\documentclass[12pt,twoside]{amsart}
\usepackage{amsmath,amssymb,mathrsfs,graphicx,bbm,amsfonts,amsthm,mathtools}
\usepackage[utf8x]{inputenc}
\usepackage{color}
\usepackage[normalem]{ulem}
\usepackage[english]{babel}
\usepackage{soul}

\usepackage{todonotes}

\usepackage{comment}

\usepackage{graphicx}
\usepackage[scriptsize,bf]{caption}

\newtheorem{theorem}{Theorem}[section]
\newtheorem{lemma}[theorem]{Lemma}
\newtheorem{proposition}[theorem]{Proposition}
\newtheorem{corollary}[theorem]{Corollary}
\newtheorem{definition}[theorem]{Definition}
\newtheorem{conjecture}[theorem]{Conjecture}
\newtheorem{remark}[theorem]{Remark}

\DeclareMathAlphabet{\mathbfit}{OML}{cmm}{b}{it}

\makeatletter
\@addtoreset{equation}{section}
\makeatother

\newcommand{\ER}{Erd\"os-Rényi\ }

\newcommand{\cadlag}{\ifmmode\mathcal D\else c\`adl\`ag \fi}

 \newcommand{\ba}{\begin{array}}
 \newcommand{\ea}{\end{array}}
 \newcommand{\bea}{\begin{eqnarray}}
 \newcommand{\eea}{\end{eqnarray}}
 \newcommand{\be}{\begin{equation}}
  \newcommand{\ee}{\end{equation}}

 \def \R {{\mathbb R}}
 \def \C {{\mathbb C}}
 \def \N {{\mathbb N}}
 \def \P {{\mathbb P}}
 
 \def \E {{\mathbb E}}

 \def \cA {\mathcal{A}}
 
 \def \cC {\mathcal{C}}
 
 \def \cE {\mathcal{E}}
 
 \def \cG {\mathcal{G}}

 \def \cM {\mathcal{M}}
 \def \cN {\mathcal{N}}
 \def \cP {\mathcal{P}}

 \def \cT {\mathcal{T}}
 \def \cU {\mathcal{U}}

 \def \a {{\alpha}}
 \def \b {{\beta}}
 
 \def \d {{\delta}}

 \def \s {{\sigma}}
 
 \def \z {{\z}}

 \def \z {{\zeta}}

 \def \G {{\Gamma}}

 \def \1{\mathbbm{1}} 

\begin{document}

\title[Limit Theorems for ERG]{Limit Theorems for Exponential Random Graphs}

\author{Alessandra Bianchi}
\address{Dipartimento di Matematica,
Universit\`a degli Studi di Padova, Via Trieste 63,
35121 Padova, Italy}\email{alessandra.bianchi@unipd.it}
\author{Francesca Collet}
\address{Dipartimento di Informatica,
Universit\`a degli Studi di Verona, Strada le Grazie~15, 
37134 Verona, Italy}\email{francesca.collet@univr.it}
\author{Elena Magnanini}
\address{WIAS, Mohrenstraße 39, 10117 Berlin, Germany}\email{elena.magnanini@wias-berlin.de}

\date{}

\begin{abstract}
We consider the edge-triangle model, a two-parameter family of exponential random graphs
in which dependence between edges is introduced through triangles.
In the so-called replica symmetric regime,
the limiting free energy exists together with a
complete characterization of the phase diagram of the model.
We borrow tools from statistical mechanics to obtain limit theorems for the edge density.
First, we investigate the asymptotic distribution of this quantity, as the graph size tends to infinity, in the various phases, and we complement this analysis with a study
of  the speed of convergence of the average edge density toward its
limiting value.
Then, we study the fluctuations of the edge density around its average value off the critical curve
 and formulate conjectures about the behavior at criticality
 based on the analysis of a mean-field approximation of the model.
Some of our results can be extended with no substantial changes to more general
classes of exponential random graphs.

    \par\bigskip\noindent
    {\it MSC 2010:} primary:  60F05, 60F15, 05C80; secondary: 60F10, 60B10.
    \par\smallskip\noindent
    {\it Keywords:} exponential random graphs, mean-field approximation, large deviations,  phase transition, standard and non-standard limit theorems, Yang-Lee theorem.
    \par\smallskip\noindent
    {\it Acknowledgements:} The authors thank Diego Alberici, Marco Formentin and Richard C. Kraaij
    for useful discussions and suggestions.
\end{abstract}

\maketitle

\section{Introduction}

In the present work, we focus on the derivation of some asymptotic properties
for the family of \emph{exponential random graphs}.
These are one of the most widely studied and promising network models
 (see \cite{RFZ} for  history), whose popularity lies in the fact
that they capture a wide variety of common network tendencies,
such as connectivity and reciprocity, by representing
a complex global structure through a set of tractable local features.
They are defined through probabilistic ensembles with one or more adjustable parameters,
and can be seen as a generalization of the classical \ER random graph \cite{ER},
obtained by adding a dependence between the random edges.
Specifically, this is realized by considering a tilted probability measure that
is proportional to the densities of certain given finite subgraphs,
in analogy to the use of potential energy
 to provide dependence between particles in a grand canonical ensemble of statistical physics.
By adjusting the specific values of these subgraph densities, through the tuning of external parameters, one can analyze the influence of various local features on the global structure of the network.

Exponential random graphs have become first popular and widely studied
in the statistical physics and network communities
(see \cite{PN3} and the surveys \cite{F1, F2}),
while in the past few years they have found space in the rigorous literature, with important contributions such as the results in \cite{CD, RY, AZ, CDey, ChNotes}.

The derivation of the free energy is a turning point for the knowledge of the model and the core of some of the aforementioned results. Indeed, it is strictly related to the cumulant generating function
of the thermodynamic observables entering in the potential energy and also a key point for obtaining large deviation results.
Its analytical expression, together with its phase diagram, is fully known
in a specific region of parameters called replica symmetric regime,
 where it can be characterized as the solution of a one-dimensional maximization problem.
Moreover, the analysis of the latter optimization problem allows to detect
a subregion in this regime where the analyticity of the free energy breaks down,
 thus revealing a phase transition \cite{CD, RY, AZ}.
As far as the asymptotics of an exponential random graph, it has been proved that,
in the replica symmetric regime, almost all realizations are close to an Erdös-Rényi random graph,
 or perhaps a finite mixture of Erdös-Rényi random graphs (phase transition) \cite{CD}.

Our analysis will be focused on the edge-triangle model (also known as Strauss's model),
a two-parameter family of exponential random graphs in which dependence between the random edges
is defined through triangles, and both edge and triangle densities are tuned by means of real parameters \cite{Str86}.
 Under the replica symmetric regime, the phase diagram for such a model is known
to contain a first order transition curve ending in a second order critical point
(qualitatively similar to the gas/liquid transition in materials) \cite{PN1, RY}.

As a first result, we  determine the asymptotic distribution of the edge density,
as the graph size $n$ tends to infinity, in the entire replica symmetric regime.
In particular, we obtain a strong law of large numbers when the parameters
are chosen outside the critical curve (see Theorem \ref{thm_LLN}),
and we prove that the edge density concentrates with high probability in a neighborhood of the maximizers of the free energy whenever working
on the critical curve (see Theorem \ref{Thm_convergence_in_distribution}).
The analysis is complemented with a study
of  the speed of convergence of the average edge density toward its
limiting value.
Specifically, we show that, up to sub-polynomial corrections,
 the speed of convergence is of order $n^{-1}$ in a subregion of parameters 
 that excludes the critical curve,
 and of order $n^{-1/2}$ at the critical point (see Proposition \ref{prop_speed}).
One main tool for the proof of these results is the large deviation principle proved
for a sequence of  \ER measures in \cite{CV11}.

We then look at the fluctuations of the edge density around its
average  for all the parameter values outside the critical curve
and off the critical point.
In particular, by exploiting properties of uniform convergence--derivable from the Yang-Lee theorem \cite{LY}--we can guarantee convergence of the moment generating function of the fluctuations and prove the central limit theorem given in Theorem \ref{thm_CLT}.

All these statements are then extended, under a proper definition of the region of criticality,
to the general family of exponential random graphs whose potential energy
is a function of different subgraph densities, including the edge density
(see Subsection \ref{generale} and theorems therein).

We conclude our work with the investigation of a simplified model
that can be seen as the mean-field approximation  of the edge-triangle model
(similar ideas can be found also in \cite{PN1}).
We show that the corresponding limiting free energy coincides with that
of the edge-triangle model, as it was also argued in \cite{MAG}.
Applying the same techniques we have adopted in the analysis of
the edge-triangle model, we then derive analogous
convergence results for the edge density of the mean-field model.
In addition, relying on its one-dimensional representation,
we are able to characterize the fluctuations of the edge density
on the critical curve (under proper conditioning) and at the critical point,
providing in this last case a non-standard central limit theorem with scaling exponent $n^{3/2}$.

These results suggest that the edge-triangle model may display
an analogous (standard versus non-standard) behavior as the parameters vary in the phase space.
In particular, supported also by the heuristics based on large deviation estimates
sketched at the end of Section \ref{speed}, we formulate  Conjecture \ref{Conj_non-standard_CLT},
about fluctuations at the critical point, and Conjectures \ref{Conj_mixture_coefficient} and
\ref{conj_conditional} about the behavior on the critical curve.

The sections are organized as follows.
In Section 2, first, we introduce the exponential random graph family and recall some main results concerning the limiting free energy. Then, we focus on the edge-triangle model; we give the definition and describe its phase diagram.
In Section 3, we state our main results on the edge-triangle model,
together with their extensions to general exponential random graphs.
Sections 4--8 are devoted to the proofs.
In Section 4, we prove properties of uniform convergence for the free energy and its derivatives; they are crucial for carrying out the proofs of Theorems \ref{thm_LLN} and \ref{thm_CLT}.
In Section 5, we provide the proof of the law of large numbers given in Theorem \ref{thm_LLN}, which is based on exponential convergence for the sequence of the edge densities.
In Section 6, we derive the central limit theorem stated in Theorem \ref{thm_CLT}
by studying the limiting behavior of the cumulant generating functions.
In Section 7, we obtain,  by means of large deviation techniques, the
concentration result, valid on the critical curve, presented in Theorem \ref{Thm_convergence_in_distribution}.
As a byproduct of this analysis, in Section 8, Proposition \ref{prop_speed},
we give bounds on the speed of convergence of the average edge density.
Section 9 is entirely devoted to the mean-field model:
we define the model, state and prove the analogs of the results
derived for the edge-triangle model in the previous sections,
sometimes in a stronger form, and then conclude with the analysis of the model
at the critical point and on the critical curve. In particular, we obtain the non-standard behavior stated in Theorem \ref{Thm_non-standard_CLT} and the conditional law of large numbers and central limit theorem stated in Theorem \ref{thm_conditional_limit_theorems_mfm} .
%

\section{Model and background}

To define the setting, let us  consider the set $\cG_n$
of all simple graphs on $n$ labeled vertices
that are identified with the elements of the set
$[n]=\{1,2,3,\ldots, n\}$.

\subsection{Exponential random graphs.}
An exponential random graph is devised to enhance or decrease the probability of specific geometric structures in the graph. Exponential weights, expressed in terms of subgraph densities, are assigned to graph ensembles. For every fixed simple graph $H$, the homomorphism density of $H$ in a graph $G$ is the probability that a random mapping $V(H) \rightarrow V(G)$, from the vertex set of $H$ to the vertex set of $G$, is edge-preserving. We write the homomorphism density as
\begin{equation}\label{def_graph_hom_density}
t(H,G) := \frac{|\text{hom}(H,G)|}{|V(G)|^{|V(H)|}}.
\end{equation}
We are going to define a probability distribution on $\mathcal{G}_n$ by means of the densities of a given selection of graphs. For any $k \in \mathbb{N}$, let $H_1, H_2, \dots, H_k$ be pre-chosen finite simple graphs (edges, stars, triangles, cycles\dots) and let $\boldsymbol{\beta}=(\beta_1, \dots, \beta_k)$ be a collection of real parameters. On $\mathcal{G}_n$ we define a real functional $\mathcal{H}_{n;\boldsymbol{\beta}}$, referred to as Hamiltonian of the model, by setting
\begin{equation}\label{Hamiltonian}
\mathcal{H}_{n;\boldsymbol{\b}}(G)=n^2\sum_{i=1}^{k}\beta_{i}t(H_{i},G)\,
\quad \mbox{ for } G\in\cG_n,
\end{equation}
and we construct the corresponding Gibbs probability density as
\be\label{eq:prob-exprg}
\mu_{n; \boldsymbol{\b}} (G)=\frac{\exp \left(\mathcal{H}_{n,\boldsymbol{\b}}(G)\right)} {Z_{n;\boldsymbol{\b}}},
\ee
where the normalizing constant $Z_{n;\boldsymbol{\beta}}$, called {\em partition function}, is given by
\be\label{eq:partit-expon}
Z_{n;\boldsymbol{\b}}=\sum_{G \in \cG_{n}} \exp \left(\mathcal{H}_{n;\boldsymbol{\b}}(G)\right) .
\ee

We will denote the related Gibbs measure and average
by $\P_{n;\boldsymbol{\b}}$ and $\E_{n;\boldsymbol{\b}}$, respectively.

\begin{remark}Notice that the Gibbs measure (exponentially)
concentrates on graphs $G\in\cG_n$ that maximize the Hamiltonian,
thus favoring or penalizing (depending on the choice of $\boldsymbol{\b}$)
the specific geometric structures entering the Hamiltonian.
\end{remark}

Finally, we define two functions that in the context
of statistical mechanics are commonly referred to as
\textit{free energies} (or pressures) of the system,
respectively, of finite and infinite size:
\be\label{FreeEnergy}
f_{n;\boldsymbol{\b}}\,:= \frac{1}{n^2}\ln Z_{n;\boldsymbol{\b}}\,\quad \text{ and } \quad
f_{\boldsymbol{\b}}\,:=\lim_{n\to+\infty} f_{n;\boldsymbol{\b}}
\,.\ee
When the function $f_{\boldsymbol{\b}}$ is well-defined,
generally depending on the value of $\boldsymbol{\b}$,
it encodes most of the asymptotic features of the model. Indeed, from its analyticity
properties, one can detect the presence
of a uniqueness/non-uniqueness phase transition for the Gibbs measure,
and in turn derive some basic connectivity properties for the large $n$ limit of the random graph.
The results given in the present section, obtained by Chatterjee and Diaconis in \cite{CD}, clarify these ideas and provide the general framework for our analysis.

Let $H_1$ be a single edge. The next theorem gives the limiting value of the free energy $f_{n;\boldsymbol{\beta}}$ whenever the parameters $\beta_2, \dots, \beta_k$ are non-negative.

\begin{theorem}[\cite{CD}, Thm.~4.1]\label{thm1_CD}
 Suppose $\beta_{2},\dots, \beta_{k}$ are non-negative. Then
\begin{equation}\label{scalar_probl}
f_{\boldsymbol{\b}} \,
=\,
 \sup_{0\leq u\leq\,1}\left(\sum_{i=1}^{k}\beta_i\,u^{E(H_i)} -\frac{1}{2}I(u)\right),
\end{equation}
where $E(H_i)$ is the number of edges in $H_i$ and $I(u):=u\ln u + (1-u)\ln(1-u)$.
 \end{theorem}

The set of values of $\boldsymbol{\b}$ such that the limiting free energy is solution of the one-dimensional variational problem \eqref{scalar_probl} is  referred to as \textbf{replica symmetric regime} (term borrowed from spin glass theory). Theorem~\ref{thm3_CD} shows that the latter phase can be slightly extended so to include (not too big) negative values for the parameters $\beta_2, \dots, \beta_k$. Theorem~\ref{thm2_CD} describes the asymptotic behavior of an exponential random graph under the replica symmetric regime.

\begin{theorem}[\cite{CD}, Thm. 4.2] \label{thm2_CD}
Suppose $\beta_{2},\dots, \beta_{k}$ are non-negative. Then, in the large $n$ limit, an exponential random graph drawn from \eqref{eq:prob-exprg} is indistinguishable from an \ER random graph with parameter $u^{*}$, where $u^{*}=u^*(\boldsymbol{\beta})$ is randomly chosen from the set of solutions of the scalar problem \eqref{scalar_probl}.
\end{theorem}

\begin{theorem}[\cite{CD}, Thm. 6.2] \label{thm3_CD}
Suppose that $\beta_{2},\dots, \beta_{k}$ are such that
\begin{equation}\label{beta_piccoli_rs}
\sum_{i=2}^{k} |\beta_{i}|E(H_{i})(E(H_{i})-1)<2,
\end{equation}
where $E(H_i)$ is the number of edges in $H_i$. Then the conclusions of  Theorems \ref{thm1_CD} and \ref{thm2_CD}
hold true. In particular, the parameters satisfying \eqref{beta_piccoli_rs} fall into the \textbf{replica symmetric regime}.
\end{theorem}

In the rest of the paper, we will restrict our analysis to the edge-triangle model, obtained when the Hamiltonian \eqref{Hamiltonian} involves only edge and triangle densities. This choice is motivated by the fact that in this setting the full characterization of the phase diagram in the replica symmetric regime is known, hence it is possible to gain control on the region of parameters where the free energy is analytic. Although all results given in Section \ref{sec_results} address specifically the edge-triangle case, some of them can be generalized by considering an arbitrary collection of subgraphs in the Hamiltonian \eqref{Hamiltonian}, as the techniques we use in the proofs are applicable provided that we  have some knowledge of the analyticity properties of the limiting free energy \eqref{scalar_probl}.
We present the statements in this general framework at the end of Section \ref{sec_results}.

\subsection{Edge-Triangle Model}
As anticipated, we now leave the general setting and we focus on the class
of exponential random graphs obtained
by considering only the contributions from edges
 and triangles. More precisely, take $G \in \mathcal{G}_n$ and fix the following subgraphs: $H_1$ is an edge and $H_2$ is a triangle. If $\beta_3=\cdots=\beta_k=0$, the Hamiltonian \eqref{Hamiltonian} reduces to
\[
\mathcal{H}_{n;\boldsymbol{\beta}} (G) = n^2 \left[ \beta_1 t(H_1,G) + \beta_2 t(H_2,G)\right].
\]
Let $E(G)$ (resp. $T(G)$) denote  the number of edges (resp.  triangles) in $G$. We have
\[
t(H_1,G) = \frac{2E(G)}{n^2} \quad \text{ and } \quad t(H_2,G) = \frac{6T(G)}{n^3}.
\]
Therefore, up to a parameter rescaling ($h=2\beta_1$; $\alpha=6\beta_2$), we can equivalently consider the Hamiltonian
 \be\label{Hamiltonian-ETmodel}
 \mathcal{H}_{n; \a,h}(G) = \frac{ \a}{ n}  T(G) +  h  E(G) \quad \text{ with } \a,h\in \R\,.
 \ee
 With a slight abuse of notation, we will denote the Gibbs probability density corresponding to \eqref{Hamiltonian-ETmodel} by $\mu_{n; \a, h}$ and the related measure and expectation by $\P_{n;\a,h}$ and $\E_{n; \a,h}$, respectively. This model is known in the literature as
edge-triangle or Strauss's model~\cite{Str86}.

\begin{remark}
Let $G=G(n,p)$ be an \ER random graph with parameters $n$ and $p$. The probability density $\mu_{n;p}^{ER}$, induced by $G$ on $\mathcal{G}_n$, is embedded in the edge-triangle model. Indeed
\be\label{eq:prob-ee}
\mu_{n;p}^{ER} (G) := (1-p)^{\binom{n}{2}} e^{ h_p E(G)} \quad \mbox{ with } \quad h_p:= \ln \frac{p}{1-p},
\ee
can be obtained from \eqref{Hamiltonian-ETmodel} by setting
{$h=h_p$} and $\a=0$.
\end{remark}

In analogy with what we have done before, if $Z_{n;\alpha,h} = \sum_{G \in \mathcal{G}_n} e^{\mathcal{H}_{n;\alpha,h}(G)}$ denotes the partition function, we write
\be\label{free_energy(a,h)}
f_{n;\a,h} := \frac{1}{n^2} \ln Z_{n;\a,h} \quad \text{ and } \quad
f_{\a,h} := \lim_{n\to+\infty}f_{n;\a,h}
\ee
for the finite and infinite volume free energies associated with \eqref{Hamiltonian-ETmodel}. We highlight that the limiting function $f_{\alpha,h}$ is well-defined only for suitable parameter values $(\alpha,h)$, as we will see in Subsection~\ref{subsect:phase_diagram}.\\
An easy computation shows that
\be
\label{averages}
\partial_h f_{n; \a,h} = \frac{\mathbb{E}_{n;\a,h}
\left(E\right)}{n^2}
\quad \text{ and } \quad
\partial_{\a} f_{n; \a,h} = \frac{\mathbb{E}_{n;\a,h}
\left(T\right)}{n^3}\,,\ee
providing a useful correspondence between the derivatives of the finite size free energy
and the averages of the edge
and of the triangle densities.

Notice that, in the present setting, condition \eqref{beta_piccoli_rs} reads $|\alpha|= 6 |\beta_2|<2$. Thus, from  Theorems \ref{thm1_CD} and \ref{thm3_CD}, it turns
out that, when $\alpha >-2$, the edge-triangle model
 is in the \textbf{replica symmetric regime}
 and therefore Theorem \ref{thm2_CD}  holds true.
 In particular, if $\alpha>-2$ the free energy exists and it is given by
\begin{equation}\label{free_energy}
f_{\a,h} \, = \,  \sup_{0 \leq u \leq 1} \left(\frac{\alpha}{6}u^3 + \frac{h}{2}u - \frac{1}{2}I(u) \right) \, = \, \frac{\a}{6}{(u^{*})}^3 +\frac{h}{2}u^{*} -\frac{1}{2}I(u^{*}),
\end{equation}
where $I(u)$ is defined in Theorem \ref{thm1_CD} and $u^{*}=u^*(\alpha,h)$ is a maximizer that solves the fixed-point equation
\begin{equation}\label{FixPointEq}
\frac{e^{\alpha\,u^{2} +h}}{1+e^{\alpha\,u^{2} +h}}=u.
\end{equation}
Thus, it is clear that, if $n$ is large, an edge-triangle model becomes indistinguishable from an \ER random graph with connection probability $u^*$. The interesting aspect is the presence of a phase transition within the replica symmetric regime. Indeed, depending on the parameters, equation \eqref{FixPointEq} can have more than one solution at which the supremum in \eqref{free_energy} is attained. Having multiplicity of optimizers translates into the possibility of having limiting \ER random graphs with very different edge densities.

The phase diagram of the edge-triangle model under the replica symmetric regime is summarized in the next subsection. The description is based on the analysis given in \cite{RY}. We would like to mention that the characterization of the (infinite size) free energy out of the replica symmetric regime is still an open question. A numerical investigation of the free energy for the case $\alpha\leq\,-2$ has been performed in \cite{GGM}.

\subsection{Phase diagram} \label{subsect:phase_diagram}
\begin{figure}[h!]
\centering
\includegraphics[scale=0.9]{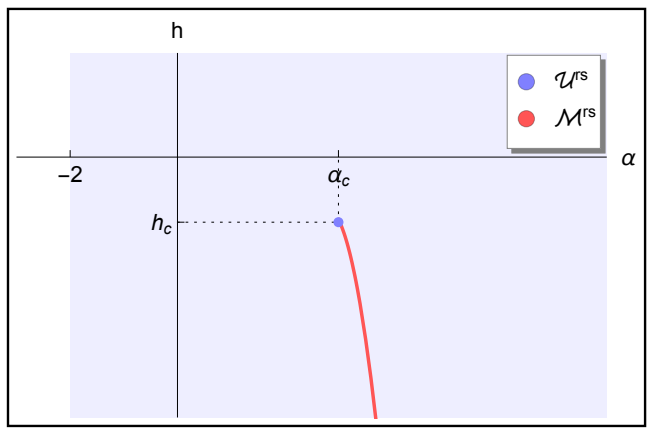}
\caption{Illustration of the phase space $(\alpha,h)$ for the edge-triangle model \eqref{Hamiltonian-ETmodel} under the replica symmetric regime. The red curve is the critical curve \eqref{curve_rs} and represents the parameter region where the optimization problem \eqref{free_energy} admits two solutions. The blue region, critical point included, corresponds to the uniqueness region for \eqref{free_energy}.
}
\label{fig:phase_diagram}
\end{figure}

For the following discussion, it is convenient to characterize the regions of the replica symmetric regime where the free energy $f_{\alpha,h}$ is well-defined, and where it is analytic (being the latter a subset of the former). It turns out that we have uniqueness for the variational problem \eqref{free_energy} on the whole replica symmetric regime except for a certain critical curve contained in the cone $\alpha > \frac{27}{8}$, $h < \ln 2 -\frac{3}{2}$. Specifically, this curve starts at the critical point $(\alpha_c,h_c) = \left(\frac{27}{8},\ln 2 -\frac{3}{2}\right)$ and it can be written as $h=q(\alpha)$ for a (non-explicit) continuous and strictly decreasing function $q$. It is known that off the curve and at the endpoint the scalar problem \eqref{free_energy} has one solution, while on the curve away from the endpoint it has two solutions  (see \cite{RY}, Prop.~3.2). To summarize, if we denote by
\begin{equation}\label{curve_rs}
\mathcal{M}^{rs} := \left\{(\alpha,h) \in (\alpha_c,+\infty) \times (-\infty,h_c): h = q(\alpha)\right\}
\end{equation}
the critical curve, then the uniqueness region within the replica symmetric regime is described by
\begin{equation}\label{reg_rs}
\mathcal{U}^{rs} := \left( (-2,+\infty) \, \times \mathbb{R} \right) \setminus \mathcal{M}^{rs}.
\end{equation}
We provide a qualitative graphical representation of the phase diagram in Fig.~\ref{fig:phase_diagram}.

The free energy $f_{\alpha,h}$ is analytic on $\mathcal{U}^{rs}\setminus \{(\alpha_c,h_c)\}$ (see \cite{RY}, Thm.~3.9). On the curve $\mathcal{M}^{rs}$ a first order phase transition occurs, and the first order partial derivatives of $f_{\alpha,h}$ have jump discontinuities. At the critical point $(\alpha_c,h_c)$ the phase transition is of the second order, and the second order partial derivatives of $f_{\alpha,h}$ diverge (see \cite{RY}, Thm.~2.1).

\section{Main results}\label{sec_results}

Let $\mathcal{E}_n$ denote the edge set of the complete graph
on $n$ vertices, with elements  labeled from 1 to $\binom{n}{2}$.
Moreover, let $\mathcal{A}_n := \{0,1\}^{\mathcal{E}_n}$. Observe that the set $\cA_n$ is equivalent to the set of the adjacency matrices of the graphs in $\cG_n$ ($n \times n$ symmetric matrices with zero diagonal entries). As a consequence, there is a one-to-one correspondence between graphs $G\in\cG_n$ and elements ${x}=(x_{i})_{i \in \cE_n}\in\cA_n$:
\begin{itemize}
\item $x_{i}=1$ if the edge $i$ is present in $G$,
\item $x_{i}=0$ otherwise.
\end{itemize}
With an abuse of nomenclature, in the sequel we will refer to the elements of $\mathcal{A}_n$ as adjacency matrices.

In view of the bijection between the sets $\mathcal{A}_n$ and $\mathcal{G}_n$, we may look at the Hamiltonian of the edge-triangle model as  a function on $\cA_n$,  defined by
\begin{equation}\label{Hamilt_ERG}
\mathcal{H}_{n;\alpha,h}(x) = \frac{\alpha}{n} \sum_{\{i,j,k\} \in \mathcal{T}_n} x_i x_j x_k + h \sum_{i \in \mathcal{E}_n} x_i,
\end{equation}
where $\mathcal{T}_n=\{\{i,j,k\} \subset \mathcal{E}_n: \{i,j,k\} \text{ is a triangle}\}$,
and then look at the corresponding Gibbs probability density
$\mu_{n;\a,h}$ as a density acting on $\cA_n$.
%
Notice that the Hamiltonian \eqref{Hamilt_ERG}
has now the form of a typical energy function used in the context of
interacting particle systems.
We then borrow some tools and techniques from statistical mechanics to analyze the model  and derive our main results. In particular, we are interested in understanding the asymptotic behavior of the number of edges.

Let $X = (X_i)_{i\in\cE_n}$ be a random element of
$\cA_n$, chosen with probability $\mu_{n;\a,h}$, and define
the (random) number of non-zero elements of $X$ as
\be\label{def:S_n}
S_n \, := \, \sum_{i\in\cE_n} X_i\,,
\ee
so that, if $x\in\cA_n$ is the adjacency matrix of the graph $G \in \mathcal{G}_n$, then $S_n(x)=E(G)$.

\begin{remark}
Notice that, even if $S_n$ is the sum of the $\binom{n}{2}$
Bernoulli random variables $X_{i}$'s, its distribution is not foregone: the $X_{i}$'s have
a non-trivial dependence structure due to the interaction Hamiltonian \eqref{Hamilt_ERG} and their distributional parameters $\E_{n;\a,h}(X_{i})$ do not have an explicit expression
as functions of the model parameters (see \cite{PN1} for further details).
\end{remark}

We consider the edge-triangle model under the replica symmetric regime and we derive classical limit theorems for the sequence $(S_n)_{n \geq 1}$. We start by characterizing the limiting distribution of the edge density. We first prove a strong law of large numbers valid for parameter values falling in the uniqueness region $\mathcal{U}^{rs}$ (see \eqref{reg_rs}). Then, we cover the phase transition curve $\mathcal{M}^{rs}$ (see \eqref{curve_rs}), showing that the edge density concentrates around the optimizers of \eqref{free_energy}.

\begin{theorem}[SLLN for $S_n$]\label{thm_LLN}
For all $(\alpha,h)\in\mathcal{U}^{rs}$, it holds
\begin{equation}\label{LLN}
\frac{2S_n}{n^2} \, \xrightarrow{\;\;\mathrm{a.s.}\;\;}{} \, u^{*}(\alpha,h) \quad \text{ w.r.t. } \mathbb{P}_{\a,h}, \text{ as } n \to +\infty,
\end{equation}	
where $u^*$ solves the maximization problem in \eqref{free_energy}.
\end{theorem}

\begin{theorem}\label{Thm_convergence_in_distribution}
For all $(\alpha,h) \in \mathcal{M}^{rs}$ and for all sufficiently small $\varepsilon>0$, there exists a constant $k=k(\varepsilon;\alpha,h)>0$ such that if
$$J(\varepsilon):= (u_1^{*}(\alpha,h)+\varepsilon,u_1^{*}(\alpha,h)-\varepsilon) \cup  (u_2^{*}(\alpha,h)+\varepsilon,u_2^{*}(\alpha,h)-\varepsilon),$$
then, for large enough $n$, it holds
\[
\mathbb{P}_{n;\alpha,h}\left(\frac{2S_n}{n^2}\in J(\varepsilon)\right)\geq 1-e^{-kn^{2}}, 
\]
where $u_1^*(\alpha,h)$ and $u_2^*(\alpha,h)$ are the two maximizers of the problem \eqref{free_energy}.
\end{theorem}

We are confident that the result presented in the previous statement can be pushed forward to a convergence in distribution to a convex combination of delta measures, but we have not been able to obtain that convergence so far.
In this respect, we present the following conjecture, which is based on  the stronger result obtainable in the mean-field setting
(see Theorem~\ref{Thm_convergence_in_distribution_mf}). 

\begin{conjecture}\label{Conj_mixture_coefficient}
For all $(\alpha,h) \in \mathcal{M}^{rs}$, it holds
\[
\frac{2S_n}{n^2} \, \xrightarrow{\;\;\mathrm{d}\;\;}{} \, \kappa \delta_{u_1^{*}}+(1-\kappa) \delta_{u_2^{*}} \quad \text{ w.r.t. } \mathbb{P}_{n;\alpha,h}, \text{ as } n \to +\infty,
\]
where $u_1^*=u_1^*(\alpha,h)$, $u_2^*=u_2^*(\alpha,h)$ solve the maximization problem in \eqref{free_energy} and the constant $\kappa \in (0,1)$ is given by
\[
\kappa=\frac{
 \sqrt{\left({1-2\alpha \left( u_1^*\right)^2(1-u_1^*)}\right)^{-1}}
}
{
 \sqrt{\left({1-2\alpha \left( u_1^*\right)^2(1-u_1^*)}\right)^{-1}}
+
 \sqrt{\left({1-2\alpha \left( u_2^*\right)^2(1-u_2^*)}\right)^{-1}}
} \,.
\]
\end{conjecture}
\vspace{0.5cm}

Having obtained a SLLN for $S_n$, it is natural to investigate at which speed the convergence occurs. The next proposition gives the convergence rate in the law of large numbers \eqref{LLN}. Our result is based on large deviation estimates and it improves the convergence rates in \cite{CDey}, Cor.~22, but it is valid for parameters lying in a smaller subregion of the uniqueness region $\mathcal{U}^{rs}$.

\begin{proposition}\label{prop_speed}
Consider $(\alpha,h) \in \mathcal{U}^{rs}$ and let $u^*$ be the unique solution of the variational problem \eqref{free_energy}. For any arbitrarily small $\varepsilon > 0$, the following assertions are true:
\begin{itemize}
\item  if $(\alpha,h) \in (-2, \alpha_c) \times \mathbb{R}$, then
\[
\frac{k_1}{n} \leq \mathbb{E}_{n;\alpha,h} \left(\left|\frac{2S_n}{n^2} - u^*(\alpha,h) \right|\right) \leq \frac{k_2}{n^{1-\varepsilon}},
\]
\item if $(\alpha,h) = (\alpha_c,h_c)$, then
\[
\frac{k_3}{n^{1/2}} \leq \mathbb{E}_{n;\alpha,h} \left(\left|\frac{2S_n}{n^2} - u^*(\alpha_c,h_c) \right|\right) \leq \frac{k_4}{n^{1/2-\varepsilon}},
\]
\end{itemize}
where $k_\cdot$'s are positive constants independent of $n$.
\end{proposition}

A second question related to the SLLN for $S_n$ concerns the behavior of the fluctuations of $S_n$ around its mean value. 
To state this type of result, we define average and variance of the edge density.

\begin{definition}
Let $X \in \mathcal{A}_n$ be an adjacency matrix randomly drawn accordingly to $\mu_{n;\alpha,h}$ and let $S_n$ be the number of non-zero elements of $X$ (see \eqref{def:S_n}).\\
For each $n \in \mathbb{N}$, we define the \textbf{average} and the \textbf{variance of the edge density}, respectively, as
\begin{equation}\label{average&variance_edge_density}
m_{n}(\a,h):=\frac{2\E_{n;\a,h}\left(S_n\right)}{n^2} \quad \text{ and } \quad v_{n}(\a,h):= \partial_h m_{n}(\a,h).
\end{equation}
\end{definition}
The reason why $v_n(\alpha,h)$ is referred to as a variance will become clear in the statement of Theorem~\ref{thm_CLT}. We can immediately notice the relevant connection between the quantities in \eqref{average&variance_edge_density} and the partial derivatives of the  free energy with respect to the parameter $h$: from \eqref{averages}, we get
\begin{equation}\label{connection_m&v_free_energy}
m_{n}(\a,h) = 2 \partial_h  f_{n;\a,h}
\quad\mbox{ and }\quad
v_{n}(\a,h)= 2 \partial_{hh} f_{n;\a,h}.
\end{equation}
	
\begin{theorem}[CLT for $S_n$]\label{thm_CLT}
For all $(\alpha,h)\in\mathcal{U}^{rs}\setminus\{(\alpha_c,h_c)\}$,
it holds
$$	
\sqrt{2} \, \frac{S_n - \frac{n^2}{2}m_{n}(\a,h)}{n} \xrightarrow{\;\;\mathrm{d}\;\;}{} \mathcal{N}(0,v(\alpha,h)) \quad \text{ w.r.t. } \mathbb{P}_{n;\a,h},\text{ as } n \to +\infty,$$
where $\mathcal{N}(0,v(\alpha,h))$ is a centered Gaussian distribution with variance given by
\begin{equation}\label{variance_CLT_ETmodel}
v(\alpha,h):=\lim_{n\to+\infty}v_{n}(\alpha,h)= \partial_h u^{*}{(\alpha,h)}.
\end{equation}
\end{theorem}

%
\begin{remark}
Thanks to a mean-field approximation of the model \eqref{Hamilt_ERG}, which will be discussed in detail in Section~\ref{sec_mean-field}, we are able to provide the explicit value of the variance in \eqref{variance_CLT_ETmodel} a posteriori. Indeed, we get $v(\alpha,h)=\frac{u^* (1-u^*)}{1-2\alpha (u^*)^2(1-u^*)}$. We refer the reader to Remark~\ref{rmk:variance_mfm} for further information.
\end{remark}

At this point, we would have liked to characterize the fluctuations of $S_n$ at the critical point $(\alpha_c,h_c)$ and on the critical curve $\mathcal{M}^{rs}$, as we are able to do so
in the mean-field setting  (see Theorems~\ref{Thm_non-standard_CLT} and~\ref{thm_conditional_limit_theorems_mfm}).
Technical difficulties have not allowed us to obtain the desired result. 
However, speculating that the edge-triangle model belongs to the same universality class of its mean-field approximation, as we believe, we are led to the following conjectures.

\begin{conjecture}[Non-standard CLT for $S_n$]\label{Conj_non-standard_CLT}
If $(\alpha,h)=(\alpha_c,h_c)$, it holds
\[	
2\, \frac{S_n - \frac{n^2}{2}  m_n(\alpha_c,h_c)}{n^{3/2}} \xrightarrow{\;\;\mathrm{d}\;\;}{} Y \quad \text{ w.r.t. } \mathbb{P}_{n;\a_c,h_c},\text{ as } n \to +\infty,
\]
where $Y$ is a random variable with Lebesgue density $\ell(y) \propto e^{-\frac{81}{64}y^{4}}$.
\end{conjecture}

We stress that the above conjecture is also supported by the heuristic computation given in \eqref{heuristic_nsCLT}, based on a large deviation property fulfilled by the law of $S_n$ under $\mathbb{P}_{n;\alpha,h}$.

To state the conjecture regarding the fluctuations of the edge-density 
on the critical curve $\mathcal{M}^{rs}$, we need some further notation.
For $(\alpha,h) \in \mathcal{M}^{rs}$, let $u_i^*(\alpha,h)$ ($i=1,2$) 
be the solutions of the scalar problem \eqref{free_energy}.
For $n\in\N$ and any fixed $\delta \in (0,1)$, consider the event
$$B_{n,\d}^{(i)}=\left\{ x\in \mathcal{A}_n \,:\,\left|\frac{2S_n(x)}{n^2} - u_i^*(\alpha,h)\right|\le n^{-\d} \right\}
$$ 
and define the conditional probability measures
\begin{equation}\label{conditional_measure}
\mathbb{P}_{n;\alpha,h}^{(i)}\left( \,\cdot\, \right) := \mathbb{P}_{n;\alpha,h} \left(\,\cdot\, \left\vert B_{n,\d}^{(i)}  \right.\right), \qquad \text{ for } i= 1,\,2\,.
\end{equation}
We denote the corresponding averages by $\mathbb{E}_{n;\alpha,h}^{(i)}$ and set  
$m_n^{(i)}(\alpha,h):=\mathbb{E}^{(i)}_{n;\alpha,h}\left(\frac{2S_n}{n^2}\right)$. 

The next conjecture provides the analog of Theorems~\ref{thm_LLN} and~\ref{thm_CLT}, but it is obtained under the constraint (\emph{conditioning}) that the edge density is close to one of the maximizers of the scalar problem \eqref{free_energy}.

\begin{conjecture}[Conditional LLN and CLT]
\label{conj_conditional}
For $i=1,2$ and for all $(\alpha,h) \in \mathcal{M}^{rs}$, it holds
\begin{equation}\label{conditional_lln}
\frac{2S_n}{n^2} \, \xrightarrow{\;\;\mathrm{a.s.}\;\;}{} \, u_i^{*}(\alpha,h)   \quad \text{ w.r.t. } \mathbb{P}^{(i)}_{n;\alpha,h}, \text{ as } n \to +\infty,
\end{equation}
and
\begin{equation}\label{conditional_clt}
\sqrt{2} \, \frac{S_n - \frac{n^2}{2} m_n^{(i)}(\alpha,h)}{n} \, \xrightarrow{\;\;\mathrm{d}\;\;}{} \, \mathcal{N}(0,v_i(\alpha,h))   \quad \text{ w.r.t. } \mathbb{P}^{(i)}_{n;\alpha,h}, \text{ as } n \to +\infty,
\end{equation}
where $\mathcal{N}(0,v_i(\alpha,h))$ is a centered Gaussian distribution with variance
\[
v_i(\alpha,h) = \frac{u_i^*(\alpha,h)[1-u_i^*(\alpha,h)]}{1-2\alpha[u_i^*(\alpha,h)]^2[1-u_i^*(\alpha,h)]}.
\]
\end{conjecture}
\vspace{0.25cm}

\subsection{Extension to the general framework}\label{generale}
The results presented in Theorems \ref{thm_LLN}, \ref{Thm_convergence_in_distribution} and \ref{thm_CLT} can be extended to the general case where the Hamiltonian \eqref{Hamiltonian} is a function of the densities of an arbitrary collection of subgraphs of the graph $G$, including the simple edge $H_1$. The main hurdle in this generality is that there is no explicit characterization of the phase transition region under the replica symmetric regime. %
%
%
However, our proofs can be mimicked (and work) if the phase diagram of the model enjoys specific characteristics. 
Let us assume that the phase structure is as follows: under the replica symmetric regime, we can identify subsets $\mathcal{U}^{rs}$ (resp. $\mathcal{M}^{rs}$), where uniqueness for the variational problem \eqref{scalar_probl} holds (resp. does not hold), and $\mathcal{C}^{rs}$, where a phase transition occurs, in the sense that the limiting free energy looses analyticity. 
In this setting, the following statements are true.

\begin{theorem}[General SLLN for $S_n$]\label{thm_gLLN}
For all $\boldsymbol{\beta}\in\mathcal{U}^{rs}$, it holds
\begin{equation}
\frac{2S_n}{n^2} \, \xrightarrow{\;\;\mathrm{a.s.}\;\;}{} \, u^{*}(\boldsymbol{\beta}) \quad \text{ w.r.t. } \mathbb{P}_{\boldsymbol{\beta}}, \text{ as } n \to +\infty,
\end{equation}	
where $u^*$ solves the maximization problem in \eqref{scalar_probl}.
\end{theorem}

\begin{theorem}\label{thm_gcritico}
Suppose that for $\boldsymbol{\beta} \in \mathcal{M}^{rs}$ the variational problem \eqref{scalar_probl} admits exactly $\ell$ solutions $u_1^*, \dots, u_{\ell}^*$. For all $\boldsymbol{\beta} \in \mathcal{M}^{rs}$ and for all sufficiently small $\varepsilon>0$, there exists a constant $k=k(\varepsilon; \boldsymbol{\beta})$ such that if
$$
J(\varepsilon):=\bigcup_{i=1}^{\ell} (u_i^{*}(\boldsymbol{\beta})+\varepsilon,u_i^{*}(\boldsymbol{\beta})-\varepsilon),
$$
then, for $n$ large enough, it holds
\[
\mathbb{P}_{n;\boldsymbol{\beta}}\left(\frac{2S_n}{n^2}\in J(\varepsilon)\right)\geq 1-e^{-kn^{2}}.
\]
\end{theorem}

\begin{theorem}[General CLT for $S_n$]\label{thm_gCLT}
For all $\boldsymbol{\beta}\in \mathcal{U}^{rs} \setminus \mathcal{C}^{rs}$,
it holds
$$	
\sqrt{2} \, \frac{S_n - \frac{n^2}{2} m_{n}(\boldsymbol{\beta})}{n} \xrightarrow{\;\;\mathrm{d}\;\;}{} \mathcal{N}(0,v(\boldsymbol{\beta})) \quad \text{ w.r.t. } \mathbb{P}_{n;\boldsymbol{\beta}},\text{ as } n \to +\infty,$$
where $\mathcal{N}(0,v(\boldsymbol{\beta}))$ is a centered Gaussian distribution with variance given by
$$v(\boldsymbol{\beta}):=\lim_{n\to+\infty}v_{n}(\boldsymbol{\beta})= \partial_{\beta_{1}} u^{*}{(\boldsymbol{\beta})},$$
being $\beta_{1}$ the first component of the vector $\boldsymbol{\beta}$.	
\end{theorem}

\section{Uniform convergence of the free energy and its derivatives}

The validity of the strong law of large numbers and of the central limit theorem as stated in Theorems \ref{thm_LLN} and \ref{thm_CLT} are deeply connected with the analyticity of the free energy and with properties of uniform convergence of its derivatives. 
Indeed, a crucial point is arguing whether it is allowed to commute limit and derivative operations. 
Typically, when dealing with spin systems, this exchange is possible due to the Griffiths, Hurst and Sherman inequality \cite{GHS}. 
Unfortunately, this result, which has been recently proved for the two-star model in \cite{BCM}, is not available for the general family of exponential random graphs. 
Our argument will rest on the Yang-Lee theorem \cite{LY}, that goes through without any difficulty in our case. 

\subsection{Polynomial representation of the partition function.}
As a first step, we consider the partition function $Z_{n;\a,h}$
on $\cA_n$ and we represent it as a polynomial. Recall that
\begin{align}
Z_{n;\a,h}&= \sum_{x\in\mathcal{A}_{n}}e^{\frac{\alpha}{n}
\sum_{\{i,j,k\}\in \cT_n}x_{i}x_{j}x_{k} + h \sum_{i\in\cE_n}x_{i}}\,.
\end{align}
To improve the readability, in the sequel
we will set $\bar n\equiv \binom{n}{2}$.
Notice that there is a bijection between $\cA_n$
and the 
power set $\cP(\cE_n)$,
that maps an element $x\in \cA_n$ to the set
$S=\{i\in\cE_n\,:\, x_{i}=1\}$.
We can then decompose $\cA_n$ in disjoint subsets as
$$\cA_n
= \bigcup_{m=0}^{\bar n} \: 
\bigcup_{S\subseteq \cE_n: |S|=m}\{x\in \cA_n\,:\, x_{i}=1 \Leftrightarrow i\in S\}, $$
and write
\begin{align}
Z_{n;\a,h}=
\sum_{m=0}^{\bar n} e^{hm} \:
\sum_{S\subseteq\cE_n: |S|=m}
e^{\frac{\alpha}{n}|\{\{i,j,k\}\subset S\,:\,
\{i,j,k\}\in\cT_n \}|}.
\end{align}
Setting $z:=e^{h}$, we arrive at the following
 representation of the partition function as a
 polynomial of degree $\bar n$:
\be\label{polinomio}
Z_{n;\a,h} \equiv Z_{\bar n}(z):=\sum_{m=0}^{\bar n}C_{m; \a}z^{m}\,.
\ee
Equivalently, letting $z_{1},z_{2},\dots, z_{\bar n}$ be the
roots of the polynomial (namely, the solutions of the equation
$Z_{\bar n}(z)=0$), we can write

\be\label{Z_radici}
Z_{\bar n}(z)=\prod_{i=1}^{\bar n}\left(1-\frac{z}{z_{i}}\right).
\ee
 We are thus in the setting analyzed
by the following Yang-Lee theorem:
 \begin{theorem}[\cite{LY}, Thm.~2]\label{teo_LY}
Let $Z_n(z)$ be the polynomial representation of
a partition function as that given in \eqref{Z_radici}.
If there exists a region $R\in\C$
 containing a segment of the real positive axis
that is always root-free then,
 as $ n\to +\infty $ and for $z\in R$,
 all quantities
\be \label{derivate}
\frac{1}{n} \ln Z_{n}(z),\quad
\frac{d^k}{d(\ln z)^k}\frac{1}{n}\ln Z_{n}(z),
\,\mbox{ with } k\in\N,
\ee
converge to analytical limits with respect to $z$.
In particular, the limit and derivative operations switch
in the whole region $R$.
\end{theorem}

The Yang-Lee theorem is a fundamental and powerful tool used in statistical mechanics to characterize the analytical properties of the infinite volume free energy relative to $Z_{n}(z)$. 
In our case, i.e. in the context of the edge-triangle model, thanks to the analysis performed in \cite{RY} (see, in particular, Thms.~2.1 and 3.9), 
we already know that the limiting free energy  $f_{\alpha,h}$ is analytic  
for all $(\alpha,h) \in \mathcal{U}^{rs} \setminus \{(\alpha_c,h_c)\}$.
As a consequence, by a contradiction argument, we can infer that the partition function \eqref{Z_radici} verifies the hypotheses of Theorem \ref{teo_LY} for all $(\alpha,h)\in \mathcal{U}^{rs}\setminus \{(\a_c,h_c)\}$. 
Hence, the convergence of the derivatives displayed in \eqref{derivate} holds.


%
\subsection{Uniform convergence of the derivatives} We show how to exploit Theorem~\ref{teo_LY} to obtain locally uniform convergence of the derivatives of the free energy.
We isolate this property in the proposition below and we briefly sketch the proof.

\begin{proposition}\label{deriv_convUnif}
Under the hypothesis of Theorem \ref{teo_LY}, the quantities displayed in \eqref{derivate} converge 
locally uniformly (in $n$) inside the region $R$.
\end{proposition}
\begin{proof}
Let $\hat\s$ be the radius of the largest open ball centered at the origin and contained in the region $R$. 
Moreover, given $\s<\hat\s$, let $C \subset R$ be the open ball centered at the origin and with radius $\sigma$. 
By construction, we have 
$|z_{i}| > \sigma$ and, for all $z \in C$, $|z|<\sigma$. 

By expanding $\frac{1}{n}\ln Z_{n}(z)$ in powers of $z$, we get
\begin{equation}
\label{serie}
\frac{1}{n}\ln Z_n(z)= \sum_{\ell=0}^{+\infty}b_{\ell}(n)z^{\ell}, 
\end{equation}
where $b_{\ell}(n)= -\frac{1}{\ell n}\sum_{i=1}^{n}\left(z_{i}\right)^{-\ell}$ 
is such that $|b_{\ell}(n)|\leq \frac{\sigma^{-\ell}}{\ell}$, for $\ell\in\N$ and uniformly 
in $n\in\N$. 
Moreover, under the hypotheses of Theorem \ref{teo_LY}, it holds
that, for all $\ell \in\N$, $b_{\ell}(\infty):= \lim_{n\to+\infty}b_{\ell}(n)$ exists (see \cite{LY}, Lem.~3) and is such that $|b_{\ell}(\infty)|\leq \frac{\sigma^{-\ell}}{\ell}$.
Putting together these observations, we obtain the following bound
\begin{align*}
\sup_{z\in C} \left| \sum_{\ell=0}^{+\infty} b_{\ell}(n) z^{\ell}
- \sum_{\ell=0}^{+\infty} b_{\ell}(\infty) z^{\ell} \right|
&\leq 
\sup_{z\in C}  \sum_{\ell=0}^{+\infty} \left| b_{\ell}(n) - b_{\ell}(\infty)\right| |z|^{\ell}\\[.2cm]
&\leq 
\sup_{z\in C}  \sum_{\ell=0}^{+\infty} \left| b_{\ell}(n) - b_{\ell}(\infty)\right| \s^{\ell}\\[.2cm]
& \leq 2\sum_{\ell=0}^{+\infty} \left(\frac{\s}{\hat\s}\right)^{\ell}\frac{1}{\ell}
\end{align*}
and the last term is finite. Taking the limit as $n$ goes to infinity, by dominated convergence, we get
\begin{align*}
\lim_{n\to+\infty}\sup_{z\in C} \left| \sum_{\ell=0}^{+\infty} b_{\ell}(n) z^{\ell}
- \sum_{\ell=0}^{+\infty} b_{\ell}(\infty) z^{\ell} \right|
&\leq 
\lim_{n\to+\infty}\sum_{\ell=0}^{+\infty} \left| b_{\ell}(n) - b_{\ell}(\infty)\right| \s^{\ell}\\[.2cm]
&= \sum_{\ell=0}^{+\infty} \lim_{n\to+\infty} \left| b_{\ell}(n) - b_{\ell}(\infty)\right| \s^{\ell}\\[.2cm]
&=0,
\end{align*}
which provides the uniform convergence of the free energy in $C$.
The very same argument can be repeated for any open ball, contained in $R$, and centered at an arbitrary point of this region. By iterating this procedure, one concludes that local uniform convergence of the free energy holds throughout the region $R$.

The local uniform convergence of the derivatives of the free energy in $R$  
can be derived analogously, as they all admit a series representation, similar to \eqref{serie}, with coefficients that are summable in $\ell \in \mathbb{N}$, uniformly in $n\in\N$.
\end{proof}

\begin{corollary}\label{cor:convergences}
Let $(\alpha,h)\in \mathcal{U}^{rs} \setminus \{(\alpha_c,h_c)\}$. If we set  $m(\alpha,h) := 2\partial_{h} f_{\alpha,h} = u^*(\alpha,h)$ and $v(\alpha,h) := 2\partial_{hh} f_{\alpha,h} = \partial_{h} u^*(\alpha,h)$, then  it holds
\[
\lim_{n \to +\infty} m_n(\alpha,h) = m(\alpha,h) \quad \text{ and } \quad \lim_{n \to +\infty} v_n(\alpha,h) = v(\alpha,h).
\]
\end{corollary}

\begin{proof}
We prove only the first display, the other being similar. The result is an immediate application of  Theorem~\ref{teo_LY} and Proposition~\ref{deriv_convUnif}, that hold true since we are working in the region $\mathcal{U}^{rs} \setminus \{(\alpha_c,h_c)\}$, where the limiting free energy is analytic (see the paragraph after Theorem \ref{teo_LY}). First notice that, since in the polynomial representation \eqref{polinomio} we have $z=e^{h}$, the derivative w.r.t. $\ln z$ of $\ln Z_{\bar n}(z)$ is equivalent, up to a constant correction, to the derivative w.r.t $h$ of the (finite size) free energy  $f_{n;\a,h}$. Therefore, the uniform convergence of the sequence $(\partial_h f_{n;\alpha,h})_{n \geq 1}$, stated in Proposition~\ref{deriv_convUnif}, allows to commute limit and derivative to get
\[
\lim_{n \to +\infty} m_n(\alpha,h) = 2\lim_{n \to +\infty} \partial_h f_{n;\alpha,h} = 2\partial_h \left[ \lim_{n \to +\infty}  f_{n;\alpha,h} \right] =  2\partial_h f_{\alpha,h} = m(\alpha,h),
\]
where the second last identity follows from Theorem~\ref{teo_LY}.
\end{proof}

\section{Strong law of large numbers via exponential convergence}

We prove the strong law of large numbers for the sequence $(S_n)_{n \geq 1}$ by exploiting exponentially fast convergence in probability. For the sake of completeness, we recall the definition of exponential convergence.

\begin{definition}[Exponential convergence]\label{def:exp_conv}
Let $(W_n)_{n \geq 1}$ be a sequence of random vectors which are defined on probability spaces $\{(\Omega_n,\mathcal{F}_n,P_n)\}_{n \geq 1}$ and which take values in $\mathbb{R}^{D}$.\\
We say that $W_n$ converges in probability exponentially fast to a constant $z_0$, and write $W_n \, \xrightarrow{\;\;\mathrm{exp}\;\;}{} \, z_{0}$, if for any $\delta>0$ there exists a number $L=L(\delta)>0$ such that
$$
P_n(\|W_n-z_{0}\|\geq\delta) \leq\,e^{-nL} \quad \text{for all sufficiently large $n$}.
$$	
\end{definition}

We present a preliminary theorem that, using the cumulant generating function, guarantees exponential convergence under general hypotheses.
\begin{theorem}[\cite{E}, Thm. II.6.3]\label{thm_Ellis}
Let $(W_n)_{n\geq 1}$ be a sequence of random vectors which are defined on probability spaces $\{(\Omega_n,\mathcal{F}_n,P_n)\}_{n\geq\,1}$ and which take values in $\mathbb{R}^{D}$. We define the cumulant generating functions as
\[
c_n(t)=\frac{1}{a_n}\ln\mathbb{E}_n[\exp(\langle t,W_n \rangle)], \qquad n=1,2,\dots, \quad t\in\mathbb{R}^{D},
\]	
where $(a_n)_{n\geq\,1}$ is a sequence of positive real numbers tending to infinity, $\mathbb{E}_n$ denotes the expectation with respect to $P_n$ and $\langle -,- \rangle$ is the Euclidean inner product on $\mathbb{R}^{D}$. We assume the following hypotheses hold:
\begin{itemize}
	\item each function $c_n(t)$ is finite for all $t\in\mathbb{R}^{D}$;
	\item $c(t)=\lim_{n\to+\infty}c_n(t)$ exists for all $t\in\mathbb{R}^{D}$ and is finite.
\end{itemize}
Then the following statements are equivalent:
\begin{itemize}
	\item $\frac{W_n}{a_n} \, \xrightarrow{\;\;\mathrm{exp}\;\;}{} \, z_{0}$;
	\item $c(t)$ is differentiable at $t=0$ and $\nabla c(0)=z_{0}$.
\end{itemize}
\end{theorem}

We rely on Theorem~\ref{thm_Ellis} to prove exponential convergence of the sequence $(S_n)_{n \geq 1}$. As we will see, existence and differentiability of the limiting cumulant generating function of the number of edges will be a direct consequence of existence and analyticity  of the infinite size free energy.
\begin{proposition}[Exponential convergence]\label{prop_exp_convergence}
For all $(\alpha,h)\in\mathcal{U}^{rs}$, it holds
\begin{equation}\label{exp_convergence}
\frac{2S_n}{n^2} \, \xrightarrow{\;\;\mathrm{exp}\;\;}{} \, u^{*}(\alpha,h) \quad \text{ w.r.t. } \mathbb{P}_{n;\a,h}, \text{ as } n \to +\infty,
\end{equation}	
where $u^*$ solves the maximization problem in \eqref{free_energy}.
\end{proposition}
\begin{proof}
Let $c_{n}(t):= 2n^{-2}\ln\mathbb{E}_{n;\a,h}[\exp(t S_n)]$ be the cumulant generating function of $S_n$ w.r.t. the Gibbs measure $\P_{n;\a,h}$. We check the requirements in Theorem~\ref{thm_Ellis}. First observe that, as $S_n$ is a bounded random variable, $c_n(t)$ is obviously finite for any $n \in \mathbb{N}$ and $t \in \mathbb{R}$. Moreover, we compute
\begin{align}\label{cgf}
c_n(t) &= \frac{2}{n^2}\ln\sum_{x\in \cA_{n}} \frac{\exp \big(\mathcal{H}_{n;\a,h}(x)+tS_n \big)}{Z_{n;\a,h}}\\[.2cm]
%
%
&=\frac{2}{n^2}\ln\frac{Z_{n;\a,h+t}}{Z_{n;\a,h}} \notag\\[.2cm]
&= 2\left(f_{n;\a,h+t}- f_{n;\a,h}\right), \notag
\end{align}
where the last equality is due to \eqref{free_energy(a,h)}.
From the existence of the infinite size free energy, given by Theorems~\ref{thm1_CD} and \ref{thm3_CD} (see also equation \eqref{free_energy}), we obtain that, for any $(\alpha, h) \in \mathcal{U}^{rs}$ and $t \in \mathbb{R}$, the limit
\begin{equation}\label{cumulant-limit}
c(t):= \lim_{n\to+\infty}c_{n}(t)= 2\left(f_{\a,h+t} - f_{\a,h}\right)
\end{equation}
exists and it is finite. Therefore, the hypotheses of Theorem \ref{thm_Ellis} are satisfied. We are left to show that $c(t)$ is differentiable at $t=0$ and that $c'(0)=u^*$, conditions equivalent to the exponential convergence statement in \eqref{exp_convergence}.
Notice that the (infinite size) free energy is analytic on $\mathcal{U}^{rs}\setminus\{(\alpha_c,h_c)\}$, and it is continuous and differentiable also at $(\alpha_c,h_c)$, since at the critical point the phase transition is of the second order (see \cite{RY}, Thm.~2.1). As a consequence, since $\partial_{t} f_{\a,h+t} = \partial_{h} f_{\a,h+t}$, we get $c'(t) = 2\partial_h f_{\alpha,h+t} = u^*(\alpha,h+t)$, for every $(\a,h+t)\in \cU^{rs}$.  In particular, we get differentiability of the function \eqref{cumulant-limit} at the origin, with $c'(0)=u^{*}(\alpha,h)$. Theorem \ref{thm_Ellis} gives the desired result.
\end{proof}
We are now ready to prove the strong law of large numbers stated in Theorem~\ref{thm_LLN}.
\begin{proof}[Proof of Theorem \ref{thm_LLN}]
Note that, as a consequence of Borel-Cantelli lemma, exponential convergence implies almost sure convergence (see \cite{E}, Thm. II.6.4). This observation together with Proposition~\ref{prop_exp_convergence} leads to the conclusion.
\end{proof}

\section{Central limit theorem}

To describe the fluctuations of $S_n$ around its mean value, we characterize the asymptotic behavior of the sequence $(V_n)_{n \geq 1}$, where
\[
V_{n} := \sqrt{2} \, \frac{S_n -\frac{n^{2}}{2}m_{n}(\alpha,h)}{n}.
\]
To this purpose, we analyze the moment generating function of the random variable $V_n$ and show that it converges to the moment generating function of a Gaussian random variable, whose variance can be explicitly computed. A key point in this strategy is to relate the moment generating function of $V_n$ to the second order derivative of the cumulant generating function of $S_n$ (recall equation \eqref{cgf}). The existence of the limit of the sequence $(c''_{n}(t))_{n\geq 1}$ for $t=t_{n}=o(1)$ will provide the variance of the limiting Gaussian.\\
We will rely on the analyticity of the free energy and on the uniform convergence of the sequence $(c''_{n}(t))_{n\geq 1}$, so Theorem \ref{teo_LY} and Proposition \ref{deriv_convUnif} will be fundamental tools in our proof. \\\\
Let $v(\alpha,h) = \partial_{h} u^*(\alpha,h)$. We want to show that the moment generating function of $V_n$ converges to the one of a Gaussian random variable with variance $v(\alpha,h)$, i.e., we want to prove that the convergence statement
\begin{equation}\label{conv:mgf}
\lim_{n\to+\infty}\E_{n;\a,h}(\exp(tV_{n})) =\exp\left(\tfrac{1}{2}v(\alpha,h)t^{2}\right)
\end{equation}
holds true for all $t\in[0,\eta)$ and some $\eta>0$. The first step is to express the average $\E_{n;\a,h}\left(\exp(tV_{n})\right)$ in terms of the second order derivative of the cumulant generating function~\eqref{cgf}. By a direct calculation, we get
\begin{equation}\label{link_c_avg&var_dens}
c'_n(t)=\frac{2\mathbb{E}_{n;\a,h+t}\left(S_{n}\right)}{n^{2}} \quad \text{ and } \quad c''_n(t)=\frac{2\text{Var}_{n;\a,h+t}(S_{n})}{n^{2}}.
\end{equation}

By using the definitions in \eqref{average&variance_edge_density}, we can connect first and second order derivatives of the cumulant generating function with average and variance of the edge density, respectively. Indeed, for all $t \in \mathbb{R}$, we have
\begin{equation}\label{link_c_avg&var_dens_2}
c'_n(t) = m_n(\alpha,h+t) \quad \text{ and } \quad c''_n(t) = v_n(\alpha,h+t).
\end{equation}

The second identity comes from the fact that, by definition of the model,
 $\partial_t m_{n}(\alpha,h+t)=
 \partial_h m_{n}(\alpha,h+t)$ and it explains why $v_n$ is referred to as a variance. In particular, we obtain $c'_n(0)=m_n(\alpha,h)$ and $c''_n(0)=v_n(\alpha,h)$. Moreover, by Corollary~\ref{cor:convergences}, we deduce
\begin{equation}\label{conv_c''_to_v}
\lim_{n \to +\infty} c''_n(0) = v(\alpha,h).
\end{equation}
Now we move back to the moment generating function of $V_n$ and we show how to write it in terms of $c''_n(t)$. Consider $t> 0$ and set $t_{n}:=\sqrt{2}t/n$. We get
\begin{align}\label{clt_decomp}
\ln\E_{n;\a,h}(\exp(tV_{n}))
%
%
&=\ln\E_{n;\a,h}\left(\exp(t_{n}S_{n})\exp\left(-\frac{t n}{\sqrt{2}} m_{n}(\alpha,h)\right)\right)\\
%
%
&= \frac{n^{2}}{2}[c_{n}(t_{n}) - t_{n}c'_{n}(0)]. \notag
\end{align}
Notice that, since $c_n(0)=0$, the term in square brackets is the difference between the function $c_n(t_n)$ and its first order Taylor expansion at zero. Therefore, by using Taylor's theorem with Lagrange remainder, one gets
\[
\ln\E_{n;\a,h}(\exp(tV_{n})) = \frac{c''_{n}(t^{*}_{n}) t^{2}}{2},
\]
for some $t^{*}_{n}\in [0,\sqrt{2}t/n]$. To conclude the proof of the central limit theorem, we need to control the limiting behavior of $c''_{n}(t^{*}_{n})$. We prove the following auxiliary lemma that leans on the uniform convergence of the sequence of derivatives
$(c'''_{n})_{n \geq 1}$.
\begin{lemma}\label{converg_dsec}
For $(\alpha,h)\in {\mathcal{U}}^{rs}\setminus \{(\alpha_c,h_c)\}$, there exists some $\eta>0$ such that $\lim_{n\to+\infty}c''_{n}(t_{n})=v(\alpha,h)$ for all $t_{n}\in[0,\eta)$ with $\lim_{n\to+\infty}t_{n} = 0$.	
\end{lemma}
\begin{proof}
Starting from the trivial bound
\[
|c''_{n}(t_{n})- v(\alpha,h)|\leq |c''_{n}(t_n) -c''_{n}(0)| + |c''_{n}(0) -v(\alpha,h)|,
\]
we observe that the second term on the right-hand side of \eqref{converg_dsec} vanishes in the limit due to \eqref{conv_c''_to_v}; while the first term can be expressed, by a Taylor expansion of first order, as $|c''_{n}(t_n) -c''_{n}(0)|=|c'''_{n}(\xi_{n})t_{n}|$, with $\xi_{n}\in[0,t_{n})$. Since $\mathcal{M}^{rs} \cup \{(\alpha_c,h_c)\}$ is a closed set in the parameter space,  it is always possible to find $\eta>0$ such that if $(\alpha,h)\in{\mathcal{U}}^{rs}\setminus\{(\alpha_c,h_c)\}$, then $(\alpha,h+t)\in{\mathcal{U}}^{rs}\setminus\{(\alpha_c,h_c)\}$ for all $t\in[0,\eta)$. As a consequence, invoking Proposition~\ref{deriv_convUnif}, we obtain that the sequence of derivatives $(c'''_{n})_{n\geq\,1}$ converges uniformly in $[0,\eta)$. In particular $c'''_{n}(\xi_{n})$ is bounded for every $n$ and consequently the term $|c'''_{n}(\xi_{n})t_{n}|$ decays to zero as $n\to+\infty$.
\end{proof}

\begin{proof}[Proof of Theorem~\ref{thm_CLT}]
Lemma~\ref{converg_dsec} readily implies that \eqref{conv:mgf} holds. The convergence in distribution of the sequence $(V_n)_{n \geq 1}$ to $\mathcal{N}(0,v(\alpha,h))$ follows from \cite{Bil86}, Sect.~30.
\end{proof}

The analysis performed in the present and in the previous section allows to extend to the critical point the convergence result for the average edge density given in Corollary~\ref{cor:convergences}.

\begin{corollary}\label{cor:convergences2}
Let $(\alpha,h)\in \mathcal{U}^{rs}$. If we set  $m(\alpha,h) := 2\partial_{h} f_{\alpha,h} = u^*(\alpha,h)$, then  it holds
\[
\lim_{n \to +\infty} m_n(\alpha,h) = m(\alpha,h).
\]
\end{corollary}
\begin{proof}
The proof is a consequence of Proposition \ref{prop_exp_convergence}
together with the following lemma:
\begin{lemma}[\cite{E}, Lem. IV.6.3]\label{conv_convex_fct}
Let $(f_{n})_{n \geq 1}$ be a sequence of convex functions on an open interval $A$ of $\R$ such that $f(t) = \lim_{n\to+\infty}f_{n}(t)$ exists for every $t\in A$.
If each $f_n$ and $f$ are differentiable at some point $t_{0}\in A$, then $\lim_{n\to+\infty}f'_{n}(t_0)$ exists
and equals $f'(t_{0})$.
\end{lemma}
The sequence $(c_n(t))_{n\geq 1}$ is a sequence of convex functions as $c_n''(t)$ can be expressed as a variance (see \eqref{link_c_avg&var_dens}).
 From the proof of Proposition \ref{prop_exp_convergence} we know that
 $c(t)= \lim_{n\to+\infty}c_{n}(t)$ exists and it is finite for any $(\alpha,h)\in \mathcal{U}^{rs}$ and $t \in \mathbb{R}$ (see \eqref{cgf} and \eqref{cumulant-limit}).
 Moreover, notice that both $c_n(t)$ and $c(t)$ are differentiable at $t=0$,
 and more precisely $c'(0)=2\partial_{h} f_{\alpha,h}= u^{*}(\alpha,h)$ and $c_n'(0)=m_n(\alpha,h)$ (see \eqref{link_c_avg&var_dens_2}).
By applying Lemma \ref{conv_convex_fct} we obtain $\lim_{n\to+\infty}c'_{n}(0)=c'(0)$,
that concludes the proof.
\end{proof}

\section{Analysis on the critical curve}
\label{Sct:critical_curve}

Our proof makes use of notions, tools, and results in the language of graph limit theory as developed in \cite{BCLSV06,BCLSV08,BCLSV12,L12,LS06}. We start the section by giving a brief overview of the concepts and properties that are relevant to our purpose.\medskip

{\em Graph limits.} Let $(G_n)_{n \geq 1}$ be a sequence of simple graphs whose number of vertices tends to infinity and let $H$ be a fixed simple graph. The sequence $(G_n)_{n \geq 1}$ converges if the graphs in the sequence become more and more similar as $n$ grows, in the sense that the homomorphism density $t(H,G_n)$ (recall definition \eqref{def_graph_hom_density}) tends to a limit $t(H)$ for every possible $H$. A reasonable limit for a sequence of graphs is thus an object from which the value $t(H)$ can be read off. Indeed, there is  a limiting object and it is not a graph, but it can rather be represented by a measurable and symmetric function $g:[0,1]^2 \rightarrow [0,1]$, called {\em graphon}. The set of all graphons is denoted by $\mathcal{W}$.\\
Any sequence of graphs that converges in the appropriate way has a graphon as limit. Conversely, every graphon arises as the limit of an appropriate graph sequence.
This limiting object determines all limits of subgraph densities: if $H$ is a simple graph with vertex set $[m]$, then
\begin{equation}\label{def_graphon_hom_density}
t(H,g) = \int_{[0,1]^m} \prod_{\{i,j\} \in \mathcal{E}(H)} g(x_i,x_j) \, dx_1 \dots dx_m,
\end{equation}
where $\mathcal{E}(H)$ denotes the edge set of $H$. A sequence of graphs $(G_n)_{n \geq 1}$ is said to converge to the graphon $g$ if, for every finite simple graph $H$,
\begin{equation}\label{density_convergence}
\lim_{n \to +\infty} t(H,G_n) = t(H,g).
\end{equation}
Intuitively, the interval $[0,1]$ represents a continuum of vertices and $g(x_i,x_j)$ corresponds to the probability of drawing the edge $\{x_i,x_j\}$. For instance, for the \ER random graph, if $p$ is fixed and $n \to +\infty$, the limit is represented by the function that is identically equal to $p$ on the unit square.
\begin{remark}\label{empiric_graphon}
Any finite simple graph admits a graphon representation. Let $H$ be a finite simple graph $H$ with vertex set $[m]$. The graphon $g^H$, corresponding to $H$, is defined by
\begin{equation}\label{graph_embedding}
g^H(x,y) = \left\{
\begin{array}{ll}
1 & \text{ if $\{\lceil mx \rceil, \lceil my \rceil\}$ is an edge in $H$}\\
0 & \text{ otherwise}
\end{array},
\right.
\end{equation}
where $(x,y) \in [0,1]^2$. In other words, $g^H$ is a step function corresponding to the adjacency matrix of $H$. This representation allows all simple graphs, regardless of the number of vertices, to be represented as elements of the single abstract space $\mathcal{W}$.
\end{remark}
In view of the above representation, the notion of convergence in terms of subgraph densities outlined above can be captured by a specific metric on $\mathcal{W}$. Let $g_1, g_2 \in \mathcal{W}$ be two graphons. Their distance is given in terms of the {\em cut distance} as
\begin{equation}\label{def_cut_distance_1}
d_{\square}(g_1,g_2) = \sup_{S,T \subseteq [0,1]} \left\vert \int_{S \times T} (g_1(x,y)-g_2(x,y)) \, dx \, dy\right\vert.
\end{equation}
A nontrivial difficulty arises from the arbitrary labeling of vertices as they are embedded in the unit interval. For this reason, we introduce an equivalence relation on $\mathcal{W}$. Let $\Sigma$ be the space of all bijections $\sigma:[0,1] \rightarrow [0,1]$ preserving the Lebesgue measure. We say that the functions $g_1,g_2 \in \mathcal{W}$ are equivalent, and we write $g_1 \sim g_2$, if $g_2(x,y) = g_1(\sigma(x),\sigma(y))$ for some $\sigma \in \Sigma$. The quotient space under $\sim$ is denoted by $\widetilde{\mathcal{W}}$ and $\tau:g \mapsto \widetilde{g}$ is the natural mapping associating a graphon with its equivalence class. Incorporating the equivalence relation $\sim$ in \eqref{def_cut_distance_1} yields a distance
\begin{equation}\label{def_cut_distance_2}
\delta_{\square} (\widetilde{g}_1,\widetilde{g}_2) = \inf_{\sigma_1, \sigma_2 \in \Sigma} \, d_{\square}(g_1(\sigma_1(x),\sigma_1(y)),g_2(\sigma_2(x),\sigma_2(y))),
\end{equation}
making $(\widetilde{\mathcal{W}},\delta_\square)$ into a compact metric space (see \cite{LS07}, Thm.~5.1).\\
To any finite graph $H$, we can associate first an element in $\mathcal{W}$: $g^H$ as in \eqref{graph_embedding}; and then an element in $\widetilde{\mathcal{W}}$: $\widetilde{g}^H = \tau (g^H)$. We get the following result for the convergence of a sequence of graphs.
\begin{theorem}[\cite{BCLSV08}, Thm.~3.8]
A sequence of graphs $(G_n)_{n \geq 1}$ converges to a limit $g \in \mathcal{W}$ in the sense of \eqref{density_convergence} if and only if $\delta_{\square}(\widetilde{g}^{G_n},\widetilde{g}) \to 0$ as $n \to +\infty$.
\end{theorem}

{\em Large deviations for random graphs.} Let $G=G(n,p)$ be an \ER random graph with parameters $n$ and $p$. The graph $G$ induces a probability distribution $\mathbb{P}_{n;p}^{ER}$ (with some abuse of notation) on $\mathcal{W}$, through the mapping $G \mapsto g^G$,  and a probability distribution $\widetilde{\mathbb{P}}_{n;p}^{ER}$ on $\widetilde{\mathcal{W}}$, through the mapping $G \mapsto g^G \mapsto \widetilde{g}^G$. The large deviation principle for the sequence of probability measures  ($\widetilde{\mathbb{P}}_{n;p}^{ER})_{n \geq 1}$ on $(\widetilde{W},\delta_{\square})$ is the main result of \cite{CV11} and it is formulated in the same way as Sanov's theorem gives a large deviation principle for a sequence of independent and identically distributed random variables.

\begin{theorem}[\cite{CV11}, Thm.~2.3]\label{Thm_LDP_ER}
For each fixed $p \in (0,1)$, the sequence $(\widetilde{\mathbb{P}}_{n;p}^{ER})_{n \geq 1}$ satisfies a large deviation principle on the space $(\widetilde{\mathcal{W}},\delta_{\square})$, with speed $n^2$ and rate function
\begin{equation}\label{ER_rate_function}
\mathcal{I}_p(\widetilde{g}) = \frac{1}{2} \int_0^1 \int_0^1  I_p(g(x,y)) \, dx \, dy,
\end{equation}
where $g$ is any representative element of the equivalence class $\widetilde{g}$ and, for $u \in [0,1]$, we set $I_p(u) = u \ln \frac{u}{p} + (1-u) \ln \frac{1-u}{1-p}$.
\end{theorem}

Building on Theorem~\ref{Thm_LDP_ER} and the following classical result, used to tilt the large deviation principle for the \ER random graph to a large deviation principle for integrals of exponential functionals w.r.t. the \ER probability distribution, we will prove Theorem~\ref{Thm_convergence_in_distribution}.

\begin{theorem}[\cite{E}, Thm.~II.7.2]\label{Thm_tilted_LDP}
Let $\mathcal{X}$ be a complete separable metric space and $(Q_n)_{n \geq 1}$ be a sequence of probability measures on the Borel $\sigma$-field $\mathscr{B}(\mathcal{X})$. Assume the sequence $(Q_n)_{n \geq 1}$ satisfies a large deviation principle with speed $a_n$ and with rate function $I$. Furthermore, let $F: \mathcal{X} \rightarrow \mathbb{R}$ be a continuous function that is bounded from above. For $n \geq 1$ and $A \in \mathscr{B}(\mathcal{X})$, define the probability measures
\[
Q_{n,F}(A) := \int_A \exp(a_nF(x)) Q_n(dx) \cdot \frac{1}{\int_{\mathcal{X}} \exp(a_n F(x)) Q_n(dx)}.
\]
Then,
\begin{enumerate}
\item
The sequence $(Q_{n,F})_{n \geq 1}$ satisfies a large deviation principle with speed $a_n$ and with rate function
\[
I_F(x) = I(x) - F(x) - \inf_{x \in \mathcal{X}} \{I(x)-F(x)\}.
\]
\item
Let $S \subset \mathcal{X}$ be a closed set that does not contain a minimizer of $I_F$. There exists a constant $k=k(S)>0$ such that
\begin{equation}\label{exp_decay_non-minimizers}
Q_{n,F}(S) \leq e^{-a_n k} \qquad \text{for all sufficiently large $n$}.
\end{equation}
\end{enumerate}
\end{theorem}

\begin{proof}[Proof \hfill of \hfill Theorem~\ref{Thm_convergence_in_distribution}]
The proof consists in showing that the sequence of the laws of the edge density w.r.t. the measure $\mathbb{P}_{n;\alpha,h}$ is exponentially tight.
In the first instance, we observe that $\mathbb{P}_{n;\alpha,h}$ can be equivalently written as a tilted probability measure on the space of graphons, that has as a priori measure the \ER one. To this aim, we need to introduce a suitable representation of the partition function. 

Consider a graph $G \in \mathcal{G}_n$ and note that,
with an abuse of notation, in view of the bijection between $\mathcal{G}_n$ and $\mathcal{A}_n$, it holds  
$\frac{2S_n(G)}{n^2} = t(H_1,G)$, being $H_1$ a single edge. 
An important property, which will be extensively used throughout the proof, is that $t(\,\cdot\,,G)= t(\,\cdot\,,\tilde{g})$, where $\widetilde{g}$ is the image in $\widetilde{\mathcal{W}}$ of the graphon representation $g^{G}$ of $G$, see Remark \eqref{empiric_graphon}. Consequently, the function $\mathcal{H}_{n;\alpha,h}$ extends naturally to $(\widetilde{W},\delta_{\square})$.
Indeed, for all $G\in\cG_n$, 
we can write 
$$
\P_{n;\alpha,h}(G)= 
\frac{\exp \left(\mathcal{H}_{n;\alpha,h}(G)\right)}{
 \, \sum_{\widetilde{g} \in \widetilde{\mathcal{W}}} \, 
 \sum_{G \in [\, \widetilde{g} \,]_n} \exp \left(\mathcal{H}_{n;\alpha,h}(G)\right)}
=
\frac{\exp \left(\mathcal{H}_{n;\alpha,h}(\widetilde{g})\right)\1(g^G \in \widetilde{g})}{
 \, \sum_{\widetilde{g} \in \widetilde{\mathcal{W}}} \, 
|[\, \widetilde{g} \,]_n|\exp \left(\mathcal{H}_{n;\alpha,h}(\widetilde{g})\right)}\,,
$$
where $[\, \widetilde{g} \,]_n:=\{G \in \mathcal{G}_n: g^G \in \widetilde{g}\}$ and $|\cdot|$ denotes the cardinality of a set. \\

Let $\mathbb{P}_{n;\frac{1}{2}}^{ER}$ be the \ER probability measure with parameters $n$ and $p=\frac{1}{2}$ and let $\widetilde{\mathbb{P}}_{n;\frac{1}{2}}^{ER}$ be the probability measure induced by $\mathbb{P}_{n;\frac{1}{2}}^{ER}$ on $\widetilde{\mathcal{W}}$. The measure $\mathbb{P}_{n;\frac{1}{2}}^{ER}$  provides a uniform probability measure on $\mathcal{G}_n$, therefore, we get
\[
\widetilde{\mathbb{P}}_{n;\frac{1}{2}}^{ER}(\widetilde{g}) = \mathbb{P}_{n;\frac{1}{2}}^{ER}([\, \widetilde{g} \,]_n) = \frac{|[\, \widetilde{g} \,]_n|}{2^{n(n-1)/2}}\,.
\]
It follows that
\begin{equation}\label{repr_pf}
\P_{n;\alpha,h}(G)=\frac{2^{-n(n-1)/2} \exp \left(\mathcal{H}_{n;\alpha,h}(\widetilde{g})\right)\1(g^G \in \widetilde{g})}
{\sum_{\widetilde{g} \in \widetilde{\mathcal{W}}} \exp\left(\mathcal{H}_{n;\alpha,h}(\widetilde{g})\right) \widetilde{\mathbb{P}}_{n;\frac{1}{2}}^{ER}(\{\widetilde{g}\})}\,.
\end{equation}

For notational convenience, let us now introduce the function
\begin{equation}
\label{U}
U_{\alpha, h}(G) = \tfrac{\alpha}{6} t(H_2,G) + \tfrac{h}{2} t(H_1,G)
\end{equation}
and write $\mathcal{H}_{n;\alpha,h}(G) = n^2 U_{\alpha, h}(G)$. For each $n \geq 1$ and each Borel set $\widetilde{A} \subseteq \widetilde{\mathcal{W}}$, we define the probability measures
\begin{equation}\label{tilted_measures}
\widetilde{\mathbb{Q}}_{n;\alpha,h}(\widetilde{A}) := \frac{\sum_{\widetilde{g} \in \widetilde{A}} \exp \left(n^2U_{\alpha, h}(\widetilde{g})\right) \widetilde{\mathbb{P}}_{n;\frac{1}{2}}^{ER}(\widetilde{g})}{\sum_{\widetilde{g} \in \widetilde{\mathcal{W}}} \exp\left(n^2U_{\alpha, h}(\widetilde{g})\right) \widetilde{\mathbb{P}}_{n;\frac{1}{2}}^{ER}(\widetilde{g})}.
\end{equation}
Since $U_{\alpha, h}$ is a continuous and bounded function on the metric space $(\widetilde{W},\delta_\square)$ (see~\cite{BCLSV08,BCLSV12}), by Theorem~\ref{Thm_tilted_LDP}(a) the sequence $\{\widetilde{\mathbb{Q}}_{n;\alpha,h}\}_{n \geq 1}$ satisfies a large deviation principle with speed $n^2$ and rate function
\begin{equation}\label{tilted_rate_function}
\mathcal{I}_{\alpha,h}(\widetilde{g}) = \mathcal{I}_{\frac{1}{2}}(\widetilde{g}) - U_{\alpha, h}(\widetilde{g}) - \inf_{\widetilde{g} \in \widetilde{\mathcal{W}}} \left\{ \mathcal{I}_{\frac{1}{2}}(\widetilde{g}) - U_{\alpha, h}(\widetilde{g}) \right\},
\end{equation}
where $\mathcal{I}_{\frac{1}{2}}$ is obtained by setting $p=\frac{1}{2}$ in \eqref{ER_rate_function}.
Since $\mathcal{I}_{\frac{1}{2}}$ is lower semicontinuous (see \cite{CV11}, Lem.~2.1), the function $\mathcal{I}_{\alpha,h}$ is too (as a sum of lower semicontinuous functions) and thus it admits a minimizer on the compact space $\widetilde{\mathcal{W}}$.
The minimizers of \eqref{tilted_rate_function} on $\widetilde{\mathcal{W}}$ are known: each optimizer is a constant function, where the constant solves the variational problem \eqref{free_energy} (see \cite{CD}, Thms. 3.1 and 4.1). In particular, whenever $(\alpha, h) \in \mathcal{M}^{rs}$, there exist two solutions $u_1^*$, $u_2^*$ for the scalar problem \eqref{free_energy}. For all sufficiently small  $\varepsilon>0$, let
\[
J(\varepsilon):= (u_1^{*}+\varepsilon,u_1^{*}-\varepsilon) \cup  (u_2^{*}+\varepsilon,u_2^{*}-\varepsilon)
\]
and consider the sets
\[
\widetilde{C}^*_{\varepsilon} := \{\tilde{g}\in \widetilde{\mathcal{W}}: t(H_1,\widetilde{g})\notin J(\varepsilon) \}
\quad  \mbox{and} \quad
C^*_{\varepsilon} :=  \left\{G\in \mathcal{G}_{n}: \frac{2S_n(G)}{n^{2}}\notin J(\varepsilon)\right\}.
\]
We point out that, due to \eqref{repr_pf} and \eqref{tilted_measures},
 it holds $\widetilde{\mathbb{Q}}_{n;\alpha,h}(\widetilde{C}^*_{\varepsilon}) = \mathbb{P}_{n;\alpha,h}\left(C^*_{\varepsilon} \right)$.
 Moreover, being $\widetilde{C}^*_{\varepsilon}$ a closed set, Theorem~\ref{Thm_tilted_LDP}(b) guarantees that, for sufficiently large $n$, there is some positive constant $k=k(\widetilde{C}^*_{\varepsilon})$ such that $\widetilde{\mathbb{Q}}_{n;\alpha,h}(\widetilde{C}^*_{\varepsilon}) \leq e^{-n^2 k}$. 
The thesis follows since 
\begin{equation}
 \mathbb{P}_{n;\alpha,h}\left(\frac{2S_n}{n^2}\in J(\varepsilon)\right)= 1- \widetilde{\mathbb{Q}}_{n;\alpha,h}(\widetilde{C}^*_{\varepsilon}) \ge 1-e^{-n^{2}k}\,.
\end{equation}
\end{proof}
The proof we did above can be carried out in the very same manner in the case when $(\alpha,h) \in \mathcal{U}^{rs}$, recovering the result in Proposition~\ref{prop_exp_convergence}. In the uniqueness regime, the set of minimizers of \eqref{tilted_rate_function} is the singleton $\widetilde{C}^*=\{\widetilde{u}^*\}$, where $\widetilde{u}^*$ is the image in $\widetilde{\mathcal{W}}$ of the unique solution $u^*$ to the scalar problem \eqref{scalar_probl}. Therefore, Theorem~\ref{Thm_tilted_LDP}(b) gives exponential convergence of the sequence $\left( 2S_n/n^2\right)_{n \geq 1}$ to $u^*$.


\begin{remark}[LDP for $\widetilde{\mathbb{Q}}_{n;\alpha,h}$]\label{Rmk_tilted_LDP}
As a byproduct of the proof of Theorem~\ref{Thm_convergence_in_distribution} we obtain that the sequence $(\widetilde{\mathbb{Q}}_{n;\alpha,h})_{n \geq 1}$ obeys a large deviation principle on the space $(\widetilde{W},\delta_{\square})$, with speed $n^2$ and rate function $\mathcal{I}_{\alpha,h}$. We have already mentioned that the rate function $\mathcal{I}_{\alpha,h}$ is lower semicontinuous; as a consequence, since $\widetilde{\mathcal{W}}$ is compact, it has compact level sets.
\end{remark}

\section{Speed of convergence}
\label{speed}

To get the speed of convergence in the law of large numbers in Theorem~\ref{thm_LLN} we find upper and lower bounds for the expectation $\mathbb{E}_{n;\alpha,h} \left( \left\vert\frac{2S_n}{n^2}-u^* \right\vert\right)$ by exploiting large deviation estimates. \\

\begin{proof}[Proof of Proposition~\ref{prop_speed}]
Take $0 < \delta < \min \{u^*, 1-u^*\}$ and introduce the event $A_\delta = \left\{ \left\vert \frac{2S_n}{n^2} - u^* \right\vert \geq \delta\right\}$. By splitting the calculation of the average over the event $A_\delta$ and its complement, we can write
\begin{equation}\label{decay_estimate}
\delta \mathbb{P}_{n;\alpha,h}(A_\delta) \leq \mathbb{E}_{n; \alpha,h} \left(\left\vert \frac{2S_n}{n^2} - u^* \right\vert\right) \leq \mathbb{P}_{n; \alpha,h} (A_\delta) + \delta.
\end{equation}
We estimate $\mathbb{P}_{n; \alpha,h} (A_\delta)$ by means of the large deviation upper and lower bounds for the sequence $(\widetilde{\mathbb{Q}}_{n;\alpha,h})_{n \geq 1}$ (see Remark~\ref{Rmk_tilted_LDP}).  Let setting and notation be as in Section~\ref{Sct:critical_curve}. Consider the closed set $\widetilde{C}_\delta = \{\widetilde{g} \in \widetilde{\mathcal{W}}: |t(H_1,\widetilde{g})-u^*| \geq \delta\} \subset \widetilde{\mathcal{W}}$ and let~$\widetilde{C}_\delta^\circ$ denote its interior. Notice that
\begin{align*}
\widetilde{\mathbb{Q}}_{n;\alpha,h}(\widetilde{C}_\delta^\circ) \leq \mathbb{P}_{n; \alpha,h}(A_\delta) &= 2^{n(n-1)/2} \, Z_{n;\alpha,h}^{-1} \sum_{\widetilde{g} \in \widetilde{C}_\delta} \exp \left( \mathcal{H}_{n; \alpha, h} (\widetilde{g})\right) \widetilde{\mathbb{P}}_{n; \frac{1}{2}}^{ER}({\widetilde{g}}) \\[.2cm]
&= \widetilde{\mathbb{Q}}_{n;\alpha,h}(\widetilde{C}_\delta),
\end{align*}
with $\widetilde{\mathbb{Q}}_{n;\alpha,h}$ the tilted probability measure defined in \eqref{tilted_measures}. The sequence $(\widetilde{\mathbb{Q}}_{n;\alpha,h})_{n \geq 1}$ satisfies a large deviation principle  with speed $n^2$ and lower semicontinuous rate function $\mathcal{I}_{\alpha,h}$, given in \eqref{tilted_rate_function}. Therefore, for sufficiently large $n$ and suitable constants $0 < k < 1 < K$, we get the following bound
\[
\exp \Big(-n^2 K \inf_{\widetilde{g} \in \widetilde{C}_\delta^\circ} \mathcal{I}_{\alpha,h}(\widetilde{g}) \Big) \leq \mathbb{P}_{n; \alpha,h}(A_\delta) \leq \exp \Big(-n^2 k \inf_{\widetilde{g} \in \widetilde{C}_\delta} \mathcal{I}_{\alpha,h}(\widetilde{g}) \Big).
\]
Since we are under the replica symmetric regime, we can characterize the solution of each variational problem in the previous display as a constant function solving the corresponding scalar problem. In other words, if \mbox{$I_{\alpha,h}: [0,1] \rightarrow \mathbb{R}$} is the function defined as $I_{\alpha,h}(x) = \frac{1}{2}I(x)  - \frac{\alpha}{6}x^3  - \frac{h}{2}x + f_{\alpha, h}$, with $I(x)=x\ln x + (1-x) \ln (1-x)$ and $f_{\alpha, h}$ given in \eqref{free_energy}, we are left with solving respectively
\[
\inf_{x \in C_\delta^\circ} I_{\alpha,h}(x) \quad \text{ and } \quad \inf_{x \in C_\delta} I_{\alpha,h}(x),
\]
where $C_\delta := \{x \in [0,1] : \vert x - u^* \vert \geq \delta\}$. The function $I_{\alpha,h}$ is continuous, positive, and convex (see \cite{RY}, Prop.~3.2); furthermore, it admits a unique zero at $x = u^*$. As a consequence, it happens
\[
\inf_{x \in C_\delta^\circ} I_{\alpha,h}(x) = \inf_{x \in C_\delta} I_{\alpha,h}(x) = \min \{I_{\alpha,h}(u^*-\delta), I_{\alpha,h}(u^*+\delta)\}>0.
\]
To conclude, we determine which is the order of the functions $I_{\alpha,h}(u^* \pm \delta)$ as $\delta$ goes to zero. We Taylor expand the functions $I_{\alpha,h}(u^* \pm \delta)$ around $u^*$ up to fourth order. Since
\[
\ln (x \pm \delta) - \ln x = \pm \frac{\delta}{x} - \frac{\delta^2}{2x^2} \pm \frac{\delta^3}{3x^3} - \frac{\delta^4}{4x^4} + o(\delta^4)\\[.3cm]
\]
and
\[
\ln \frac{x \pm \delta}{1-x \mp \delta} = \ln \frac{x}{1-x} \pm \frac{\delta}{x(1-x)} + \frac{(2x-1)\delta^2}{2x^2(1-x)^2} \pm  \frac{(3x^2-3x+1)\delta^3}{3x^3(1-x)^3} + o(\delta^3),
\]
we obtain
\begin{align}\label{rate_fct_expansion}
I_{\alpha,h}(u^* \pm \delta) &= \pm \frac{\delta}{2} \left[ \ln \frac{u^*}{1-u^*} - \alpha (u^*)^2 - h\right] + \frac{\delta^2}{2} \left[ \frac{1}{2u^*(1-u^*)} - \alpha u^*\right] \\[.3cm]
&\quad \pm \frac{\delta^3}{2} \left[ \frac{2u^*-1}{6(u^*)^2(1-u^*)^2} - \frac{\alpha}{3}\right] + \frac{\delta^4}{24} \, \frac{3(u^*)^2 -3u^*+1}{(u^*)^3(1-u^*)^3} + o(\delta^4).\nonumber
\end{align}
Notice that the coefficient of the first order term in \eqref{rate_fct_expansion} vanishes for all values of the parameters in the replica symmetric regime, as $u^*$ is solution to the fixed point equation \eqref{FixPointEq}. Moreover,
\begin{itemize}
\item if $(\alpha,h) \in (-2, \alpha_c) \times \mathbb{R}$, then the first non-vanishing coefficient in \eqref{rate_fct_expansion} is the coefficient of the second order term, which equals $\frac{I_{\alpha,h}''(u^*)}{2}$ and it is strictly positive as $I_{\alpha,h}$ is strictly convex in the parameter range we are considering. Therefore, we get $I_{\alpha,h}(u^* \pm \delta) = \frac{I_{\alpha,h}''(u^*)}{2} \, \delta^2 + o(\delta^2)$ for small $\delta$.
\vspace{.2cm}
\item if $(\alpha,h) = (\alpha_c,h_c) = \left( \frac{27}{8}, \ln 2 - \frac{3}{2} \right)$, since $u^*=u^*(\alpha_c,h_c) = \frac{2}{3}$, the first non-vanishing coefficient in \eqref{rate_fct_expansion} is the coefficient of the fourth order term, which gives $\frac{81}{64}$. Hence, we get $I_{\alpha,h}(u^* \pm \delta) = \frac{81}{64} \, \delta^4 + o(\delta^4)$ for small $\delta$.
\end{itemize}
Going back to \eqref{decay_estimate}, we can develop further the estimate and write
\begin{equation}\label{decay_estimate_final}
\delta \exp(-n^2KL(\delta)) \leq \mathbb{E}_{n; \alpha,h} \left(\left\vert \frac{2S_n}{n^2} - u^* \right\vert\right) \leq \exp(-n^2kL(\delta)) + \delta,
\end{equation}
with $L(\delta)$ a positive constant. Moreover, the analysis we have just carried out  shows that, for small $\delta$, $L(\delta) \in O(\delta^2)$ if $(\alpha,h) \in (-2,\alpha_c) \times \mathbb{R}$, whereas $L(\delta) \in O(\delta^4)$ if $(\alpha,h) = (\alpha_c,h_c)$. In particular, by choosing $\delta=\delta(n)=n^{-1}$ and $\delta=\delta(n)=n^{-\gamma}$, with $0 < \gamma < 1$, we obtain the tightest possible lower and upper bounds in non-critical cases; whereas, the choices  $\delta=\delta(n)=n^{-1/2}$ and $\delta=\delta(n)=n^{-\gamma}$, with $0 < \gamma < \frac{1}{2}$, give the optimal bounds at criticality.
\end{proof}

\begin{remark}[Heuristics on non-standard critical behavior]
The Taylor expansion \eqref{rate_fct_expansion} of the function $I_{\alpha,h}$ around its minimum, explicitly evaluated at the critical point $(\alpha_c,h_c)$,
allow performing a heuristic calculation to support Conjecture~\ref{Conj_non-standard_CLT}. 
We define the random variable
\[
U_n := 2 \, \frac{S_n-\frac{n^2}{2}u^*(\alpha_c,h_c)}{n^{3/2}},
\]
so that the target random variable in the conjecture can be written as
\be\label{shift}
2 \, \frac{S_n-\frac{n^2}{2}m_n(\alpha_c,h_c)}{n^{3/2}}=
U_n + \sqrt{n}(u^*-m_n(\alpha_c,h_c))\,.
\ee

We use the approximation $\mathbb{P}_{n;\alpha,h} \left( \frac{2S_n}{n^2} \approx u \right) = \exp  \left( -n^2 I_{\alpha,h}(u) + o(n^2) \right)$, whenever $u$ is closed to $u^*$, and we Taylor expand $I_{\alpha,h}$ around $u^*$ up to fourth order. As $n \to +\infty$, we find
\begin{align}\label{heuristic_nsCLT}
\mathbb{P}_{n;\alpha_c,h_c} \left( U_n \in dx \right) &= \mathbb{P}_{n;\alpha_c,h_c} \left( \tfrac{2S_n}{n^2} \in u^*+\tfrac{dx}{\sqrt{n}} \right) \\[.2cm]
&\approx e^{ -n^2 I_{\alpha_c,h_c} \left( u^* + \frac{x}{\sqrt{n}}\right) + o(n^2)} dx = e^{-\frac{81}{64} x^4 + o(n^2)}dx, \nonumber
\end{align}
which is the same density $\ell^c$ found in Theorem~\ref{Thm_non-standard_CLT}.
However, notice that the error term appearing in \eqref{heuristic_nsCLT}
might be relevant, also accordingly to Theorem 1.4(c) in \cite{MX} 
stated in the context of the two-star model,  as well as the shift term in \eqref{shift} (see Proposition \ref{prop_speed}). 
We believe that a subtle compensation among this two contributions produces the conjectured result.
\end{remark}

\section{Mean-field approximation}
\label{sec_mean-field}

In this section, we consider a mean-field approximation of the edge-triangle model \eqref{Hamilt_ERG}. The mean-field approximation provides one of the simplest and most common schemes for analyzing a model, when exact computations are infeasible, and it is a helpful tool for deriving heuristic results. It will be shown that all results valid for the distribution $\mathbb{P}_{n;\alpha,h}$, are still valid in this approximation even in a stronger form (e.g. Theorem \ref{Thm_convergence_in_distribution_mf} and Proposition \ref{prop_speed_mmf}).
Moreover, we will derive a non-standard CLT at the critical point $(\alpha_c,h_c)$,
that together with the heuristics at the end of the previous section
suggests Conjecture \ref{Conj_non-standard_CLT}. 
We will then focus on the behavior of the model on the critical curve,
and prove conditional LLN and CLT for the edge density.
We will conclude the section with a comparison among the edge-triangle model and its mean-field approximation.

\subsection{Mean-field model and limiting free energy}
For any given $\alpha,h\in\mathbb{R}$, let us
consider the \emph{mean-field Hamiltonian}
\begin{equation}\label{H_mf}
\bar{\mathcal{H}}_{n;\alpha,h}(x):=\frac{4\alpha}{3n^4}\left(\sum_{i\in\cE_n} x_{i}\right)^{3}+h\sum_{i\in\cE_n} x_{i}\,, \quad \text{ for } x\in\mathcal{A}_n \,,
\end{equation}
with associated Gibbs probability density
\begin{equation}\label{mean-field-mes}
\bar{\mu}_{n;\alpha,h}(x):=
\frac{e^{\bar{\mathcal{H}}_{n;\alpha,h}(x)}}{\bar{Z}_{n;\alpha,h}}\,,
\qquad \bar{Z}_{n;\alpha,h}:=\sum_{x\in\cA_n}e^{\bar{\mathcal{H}}_{n;\alpha,h}(x)}.
\end{equation}
We denote the corresponding measure and expectation
by $\bar{\P}_{n;\a,h}$ and $\bar{\E}_{n;\a,h}$,
 respectively. Moreover, as usual, we define the finite size free energy as
\begin{equation}\label{fe_mf}
\bar{f}_{n;\alpha,h}:= \frac{1}{n^{2}}\ln \bar{Z}_{n;\alpha,h}.
\end{equation}

Notice that this is a mean-field model in the sense
that the Hamiltonian (and the corresponding probability density $\bar{\mu}_{n;\alpha,h}$)
can be seen as a function of the one dimensional
parameter
$\frac{2}{n^2}S_n(x) = \frac{2}{n^2}\sum_{i\in\cE_n}x_{i}$,
taking values in $\G_{n} :=\left\{0,\frac{2}{n^2}, \dots,1-\frac{1}{n}\right\}$. In particular, for all $x\in\cA_n$ such that $\frac{2S_n(x)}{n^2}= m$, with $m\in\G_n$,
we have
\be\label{H_mf2}
\bar{\mathcal{H}}_{n;\alpha,h}(x)= \bar{\mathcal{H}}_{n;\alpha,h}(m)= n^2 \left(\frac{\alpha}{6} m^{3} + \frac{h}{2} m \right).
\ee
This identity induces on $\G_n$ the following measure
\be\label{Gibbs_m}
\bar{\P}_{n;\a,h}\left(\frac{2S_n}{n^2}\in A\right)
= \sum_{m\in A } \cN_m \frac{e^{n^2 \left( \frac{\alpha}{6}m^3+\frac{h}{2}m\right)}}{\bar{Z}_{n;\a,h}}\,,\qquad \text{ for } A\subseteq \G_n\,,
\ee
where $\mathcal{N}_{m}:= \binom{\frac{n(n-1)}{2}}{\frac{n(n-1)m}{2}}$ corresponds to the number of graphs $x\in\cA_n$ with edge density $2S_n(x)/n^2 =m$.

This model provides a mean-field approximation of the edge-triangle model \eqref{Hamilt_ERG} in the following sense: whenever the parameters $\alpha$, $h$ are chosen so that the edge-triangle model is in the replica symmetric regime, the two models share the same limiting free energy. We will show this result below in Theorem~\ref{fe_meanf_exp}. Before we prove an auxiliary lemma that will be a useful tool in the proof of Theorem~\ref{fe_meanf_exp} and many others in this section. It says that the main contribution to $\bar{Z}_{n;\alpha,h}$ comes from the sum over the neighborhood(s) of the solution(s) of the scalar problem \eqref{free_energy}.\\
The definitions of regions $\mathcal{U}^{rs}$ and $\mathcal{M}^{rs}$ in the parameter space are given in Subsection~\ref{subsect:phase_diagram}.

\begin{lemma}\label{lemma_Z}
Let $(\alpha,h) \in (-2,+\infty) \times \mathbb{R}$ and let $f_{\alpha,h}$ be the infinite volume free energy  of the edge-triangle model given in \eqref{free_energy}. Let $g_{\a,h}:[0,1] \rightarrow \mathbb{R}$ be the function defined as
\be\label{g}
g_{\a,h}(m):= \frac{\alpha}{6}m^{3} + \frac{h}{2}m - \frac{I(m)}{2}.
\ee
Finally, let $\d\in \left(0,\frac{3}{8} \right)$ be fixed and set $R^{(n)}:=\left\{ -n^{1-\d}, -n^{1-\d}+\frac{2}{n},\ldots, n^{1-\d} \right\}$ and $R_c^{(n)}:= \left\{-n^{1/2-\delta}, -n^{1/2-\delta}+\frac{2}{n^{3/2}}, \dots, n^{1/2-\delta}\right\}$. It holds true that, as $n \to +\infty$,
\be\label{claim}
\bar{Z}_{n;\a,h}=  \frac{e^{n^2 f_{\a,h}}}{n\sqrt{\pi }} \left( D^{(n)}(\alpha,h)\right) (1+o(1)),
\ee
where
\begin{itemize}
\item if $(\alpha,h) \in \mathcal{U}^{rs}\setminus \{(\a_c,h_c)\}$,
\[
D^{(n)}(\alpha,h) := D_*^{(n)} = \sum_{x\in R^{(n)}} \frac{e^{- c_0 x^2 + \frac{k_0}{n} x^3}}
{\sqrt{\left(u^*+\frac{x}{n}\right) \left(1-u^* -\frac{x}{n}\right)}},
\]
with $u^*$ solution of the scalar problem \eqref{free_energy} and the constants defined as
\[
c_0 := -\frac{g''_{\a,h} (u^*)}{2} = \frac{1-2\alpha(u^*)^2(1-u^*)}{4u^*(1-u^*)} >0
\]
and $k_0:= g'''_{\a,h} (\tilde{u})/6$, for some $\tilde{u}$ such that $|\tilde{u}-u^*|< n^{-\d}$;\\
\item if $(\alpha,h) \in \mathcal{M}^{rs}$,
\[
D^{(n)}(\alpha,h) := D_1^{(n)} + D_2^{(n)},
\]
with
\begin{align*}
D_1^{(n)} &:= \sum_{x\in R^{(n)}} \frac{e^{- c_1 x^2 + \frac{k_1}{n} x^3}}
{\sqrt{\left(u_1^*+\frac{x}{n}\right) \left(1-u_1^* -\frac{x}{n}\right)}}\\[.2cm]
D_2^{(n)}&:= \sum_{x\in R^{(n)}} \frac{e^{- c_2 x^2 + \frac{k_2}{n} x^3}}
{\sqrt{\left(u_2^*+\frac{x}{n}\right) \left(1-u_2^* -\frac{x}{n}\right)}},
\end{align*}
where $u_1^*$, $u_2^*$ are solutions of the scalar problem \eqref{free_energy} and the constants are defined as 
\[
c_i := - g''_{\a,h} (u_i^*)/2 = \frac{1-2\alpha(u_i^*)^2(1-u_i^*)}{4u_i^*(1-u_i^*)} >0
\]
and $k_i:= g'''_{\a,h} (\tilde{u})/6$, for some $\tilde{u}$ such that $|\tilde{u}-u_i^*|< n^{-\d}$ ($i=1, 2$);\\
\item if $(\alpha,h) = (\alpha_c, h_c)$,
\[
D^{(n)}(\alpha_c,h_c) := D_c^{(n)} = \sum_{x\in R_c^{(n)}} \frac{e^{- \frac{81}{64} x^4 + \frac{k_c}{\sqrt{n}} x^5}}
{\sqrt{\left(u^*+\frac{x}{\sqrt{n}}\right) \left(1-u^* -\frac{x}{\sqrt{n}}\right)}},
\]
with $u^*=u^*(\alpha_c,h_c)=\frac{2}{3}$ solution of the scalar problem \eqref{free_energy} and $k_c:=g^{(v)}_{\a_c,h_c} (\tilde{u})/5!$, for some $\tilde{u}$ such that $|\tilde{u}-u^*|< n^{-\d}$.
\end{itemize}
\end{lemma}

\begin{proof}
In view of the representation \eqref{H_mf2}, we write the partition function as
\begin{equation}\label{Z_mf}
\bar{Z}_{n;\alpha,h}=  \sum_{m\in\G_{n}} \mathcal{N}_{m}
e^{n^2 (\frac{\alpha}{6}m^{3} + \frac{h}{2}m)}\,.
\end{equation}
The coefficient $\cN_m$ can be
approximated by using Stirling's formula.
In particular, there exist two constants $ c,C>0$
such that, for all $m\in\G_n$, it holds
\begin{equation}\label{Stirling-rough}
cn^{-1} e^{-\frac{n^2}{2}I(m)}
\leq \mathcal{N}_{m}\leq C n e^{-\frac{n^2}{2}I(m)}\,.
\end{equation}
The estimate \eqref{Stirling-rough} can  be made more precise assuming
$ n^{-2}\ll m \ll 1-n^{-2}$, so that $n^2 m\to+\infty$ and $n^2(1-m)\to+\infty$.
We thus obtain
\begin{equation}\label{Stirling}
\mathcal{N}_{m}= \frac{e^{-\frac{n^2}{2}I(m)}}{n \sqrt{\pi m(1-m)}} (1+o(1))\,.
\end{equation}
We start by considering $(\a,h)\in\cM^{rs}$.
Let $u_1^*$ and $u_2^*$ be the solutions of the scalar problem \eqref{free_energy}.
For any large enough $n$, we can split the sum in \eqref{Z_mf}
in three parts, accordingly to the partition
$\G_n= B_{u_1^*} \cup B_{u_2^*} \cup \,\cC$, where
\be\label{neighborhoods}
\begin{split}
&B_{u_1^*} \equiv B_{u_1^*}(n,\d) := \{m\in\G_n\,:\, |m-u_1^*|\leq n^{-\d}\},\\
&B_{u_2^*} \equiv B_{u_2^*}(n,\d) := \{m\in\G_n\,:\, |m-u_2^*|\leq n^{-\d}\},\\
&\cC \equiv \mathcal{C}(n,\d) := \G_n \setminus \left(B_{u_1^*}\cup B_{u_2^*}\right),
\end{split}
\ee
and consider the three sums separately.
\vspace{0.5cm}

\noindent \emph{Neighborhood of $u_1^*$.}
From the Stirling's approximation \eqref{Stirling}, we first get
\be\label{somma0}
\bar{Z}_{n;\a,h}(B_{u_1^*}):=\sum_{m\in B_{u_1^*}}\cN_m e^{n^2(\frac{\alpha}{6}m^{3} + \frac{h}{2}m)}
= \frac{1}{n\sqrt{\pi}}
\sum_{m\in B_{u_1^*}}\frac{e^{n^2g_{\a,h}(m)}}{\sqrt{m(1-m)}} (1+o(1)).
\ee
We Taylor expand $g_{\a,h}(m)$ around $u_1^*$.
Since  $g'_{\a,h}(u_1^{*})=0$,
as $u_1^{*}$ is a maximizer of $g_{\a,h}$,
it holds that
\be\label{taylor_g_M}
g_{\a,h}(m)=
g_{\a,h}(u_1^*)-c_1(m-u_1^*)^2 + k_1(m-u_1^*)^3,
\ee
where $c_1>0$ and $k_1$ are the constants given in the statement.
Inserting \eqref{taylor_g_M} in \eqref{somma0}, and recalling that
$ g_{\a,h}(u_1^*)\equiv f_{\a,h}$, yields
\be\label{somma1}
\begin{split}
\bar{Z}_{n;\a,h}(B_{u_1^*})
&= \frac{e^{n^2 f_{\a,h}}}{n\sqrt{\pi}}
\sum_{m\in B_{u_1^*}}\frac{e^{-c_1 n^2(m-u_1^*)^2 + k_1 n^2(m-u_1^*)^3}}{\sqrt{m(1-m)}} (1+o(1))\\
& = \frac{e^{n^2 f_{\a,h}}}{n\sqrt{\pi}} D_1^{(n)} (1+o(1)),
\end{split}
\ee
where the last identity follows from the change of variable $x=n(m-u_1^*)$.
\vspace{0.5cm}

\noindent\emph{Neighborhood of $u_2^*$.}
Since $u_2^*$ is the maximizer of $g_{\a,h}$ over $B_{u_2^*}$,
following the same procedure as for the sum over $B_{u_1^*}$, we get that
\be\label{somma2}
\begin{split}
\bar{Z}_{n;\a,h}(B_{u_2^*})
&= \frac{e^{n^2 f_{\a,h}}}{n\sqrt{\pi}}
\sum_{m\in B_{u_2^*}}\frac{e^{-c_2 n^2(m-u_1^*)^2 + k_2 n^2(m-u_1^*)^3}}{\sqrt{m(1-m)}} (1+o(1))\\
& = \frac{e^{n^2 f_{\a,h}}}{n\sqrt{\pi}} D_2^{(n)} (1+o(1)),
\end{split}
\ee
with constants $c_2>0$ and $k_2$ as given in the statement.

Notice that since $D_1^{(n)}$ and $ D_2^{(n)}$ are positive and uniformly bounded
in $n\in\N$,  it turns out  that
$\bar{Z}_{n;\a,h}(B_{u_1^*})$ and $\bar{Z}_{n;\a,h}(B_{u_2^*})$ have the same order
as $n\to+\infty$.
\vspace{0.5cm}

\noindent\emph{Complement set $\cC$.}
From the Stirling's approximation \eqref{Stirling-rough},  we get
\be\label{somma3}
\begin{split}
\bar{Z}_{n;\a,h}(\cC)
&\leq C n \sum_{m\in \cC} e^{n^2 g_{\a,h}(m)}\leq
\frac{C}{2} n^3 e^{n^2 f_{\a,h}} e^{-n^2( f_{\a,h} -\max_{m\in\cC} g_{\a,h}(m))}\\
& \leq e^{n^2 f_{\a,h}} e^{-k n^{2-2\d}}= e^{n^2 f_{\a,h}}o(1),
\end{split}
\ee
where the last inequality holds for some $k>0$ and any  large enough $n$.
We used  that $|m- u_j^*|> n^{-\d}$, for $j=1,2$ and for all $m\in\cC$,
and that $2-2\d>0$, whenever $\d<1$.

Taking  \eqref{somma1}-\eqref{somma3} together gives
$$
\bar{Z}_{n;\a,h}= \left(\bar{Z}_{n;\a,h}(B_{u_1^*})+\bar{Z}_{n;\a,h}(B_{u_2^*})\right)
(1+o(1))=  \frac{e^{n^2 f_{\a,h}}}{n\sqrt{\pi}}
\left(D_1^{(n)}+D_2^{(n)}\right) (1+o(1)),
$$
that concludes the proof of the statement for any $(\a,h)\in\cM^{rs}$.

If $(\a,h)\in\cU^{rs}$,
 being $u^*$ the unique solution of the scalar problem \eqref{free_energy},
we can proceed similarly by splitting the sum in \eqref{Z_mf}
over the neighborhood of $u^*$ and its complement.
The  previous argument then applies straightforwardly
by just adapting the Taylor's approximation
of $g_{\a,h}(m)$ around its maximum.
By direct computation, and taking into account that
$\alpha_{c}=\frac{27}{8}$, $h_{c}= \ln 2 -\frac{3}{2}$ and $u^{*}(\alpha_c,h_{c})=\frac{2}{3}$,
we get
\be\label{taylor_g}
g_{\a,h}(m)-g_{\a,h}(u^*)=
\left\{
\begin{array}{ll}
-c_0(m-u^*)^2 + k_0(m-u^*)^3 &
\mbox{if } (\alpha,h) \in \mathcal{U}^{rs}\setminus \{(\a_c,h_c)\}\\
\\
-\frac{81}{64}(m-u^*)^4 + k_c (m-u^*)^5 &
 \mbox{if } (\alpha,h) = (\a_c,h_c),
\end{array}
\right.
\ee
where the constants $c_0 > 0$, $k_0$, and $k_c$ are those
given in the statement.
As a consequence, after the change of variable
$x=n(m-u^*)$ if  $(\a,h)\in\cU^{rs}\setminus\{(\a_c,h_c)\}$,
or  \mbox{$x=\sqrt{n}(m-u^*)$} if $(\a,h)=(\a_c,h_c)$, we get
$$
\bar{Z}_{n;\a,h}=
\bar{Z}_{n;\a,h}(B_{u^*})(1+o(1))=
\left\{
\begin{array}{ll}
  \frac{e^{n^2 f_{\a,h}}}{n\sqrt{\pi}}
D_*^{(n)} (1+o(1)) & \mbox{ if } (\a,h)\in\cU^{rs}\setminus\{(\a_c,h_c)\}
\\[.4cm]
 \frac{e^{n^2 f_{\a,h}}}{n\sqrt{\pi}}
D_c^{(n)} (1+o(1)) & \mbox{ if } (\a,h)=(\a_c,h_c),
\end{array}\right.$$
that concludes the proof.
We point out that, when at the critical point, the definition of $B_{u^*}$ as $n^{-\delta}$-neighborhood of $u^*$, under the assumption $ 0 < \delta < \frac{3}{8}$, guarantees that the contribution to the partition function of the term $\bar{Z}_{n;\alpha_c,h_c}(\mathcal{C})$ is exponentially small in $n$, thus negligible.

\end{proof}

\begin{theorem}\label{fe_meanf_exp}
Let $(\alpha,h) \in (-2,+\infty) \times \mathbb{R}$ and let $f_{\alpha,h}$ be the infinite volume free energy  of the edge-triangle model given in \eqref{free_energy}.
Then it holds
\begin{equation}
\label{equiv_fren}
\lim_{n\to+\infty}\bar{f}_{n;\alpha,h}= f_{\alpha,h}\,.
\end{equation}
\end{theorem}
This result was already observed in \cite{MAG}, Ch.~IV, and it is a straightforward consequence of \eqref{claim}.

\subsection{Limit theorems for the mean-field model}
All results presented in Section \ref{sec_results} can be easily transferred, and in some cases strengthened, to the mean-field setting. We summarize the statements and the relevant proofs below.

\begin{theorem}[SLLN for $S_n$ w.r.t. $\bar{\mathbb{P}}_{n;\alpha,h}$]\label{thm_LLN_mf}
For all $(\alpha,h)\in\mathcal{U}^{rs}$, it holds
\[
\frac{2S_n}{n^2} \, \xrightarrow{\;\;\mathrm{a.s.}\;\;}{} \, u^{*}(\alpha,h) \quad \text{ w.r.t. } \bar{\mathbb{P}}_{\alpha,h}, \text{ as } n \to +\infty,
\]
where $u^*$ solves the maximization problem in \eqref{free_energy}.
\end{theorem}

\begin{theorem}\label{Thm_convergence_in_distribution_mf}
For all $(\alpha,h) \in \mathcal{M}^{rs}$, it holds
\[
\frac{2S_n}{n^2} \, \xrightarrow{\;\;\mathrm{d}\;\;}{} \, \kappa \delta_{u_1^{*}(\alpha,h)}+
(1-\kappa)\delta_{u_2^{*}(\alpha,h)}  \quad \text{ w.r.t. } \bar{\mathbb{P}}_{n;\alpha,h}, \text{ as } n \to +\infty,
\]
where $u_1^*$, $u_2^*$ solve the maximization problem in \eqref{free_energy}, and
$$\kappa=\frac{
 \sqrt{\left({1-2\alpha \left( u_1^*\right)^2(1-u_1^*)}\right)^{-1}}
}
{
 \sqrt{\left({1-2\alpha \left( u_1^*\right)^2(1-u_1^*)}\right)^{-1}}
+
 \sqrt{\left({1-2\alpha \left( u_2^*\right)^2(1-u_2^*)}\right)^{-1}}
} \,.
$$
\end{theorem}

\begin{theorem}[CLT for $S_n$ w.r.t. $\bar{\mathbb{P}}_{n;\alpha,h}$]\label{thm_CLT_mf}
For all $(\alpha,h)\in\mathcal{U}^{rs}\setminus\{(\alpha_c,h_c)\}$,
it holds
\begin{equation}\label{CLT_mf}
\sqrt{2} \, \frac{S_n - \frac{n^2}{2}\bar{m}_{n}(\a,h)}{n} \xrightarrow{\;\;\mathrm{d}\;\;}{} \mathcal{N}(0,v(\alpha,h)) \quad \text{ w.r.t. } \bar{\mathbb{P}}_{n;\alpha,h},\text{ as } n \to +\infty,
\end{equation}
where $\bar{m}_{n}:=\bar{\mathbb{E}}_{n;\alpha,h}(\frac{2S_n}{n^2})$
 and $\mathcal{N}(0,v(\alpha,h))$ is a  centered Gaussian distribution with variance
\begin{equation}\label{variance_CLT_mfm}
v(\alpha,h):=\lim_{n\to+\infty}v_{n}(\alpha,h)= \partial_h u^{*}{(\alpha,h)}.
\end{equation}
\end{theorem}

The proofs of Theorems~\ref{thm_LLN_mf} and~\ref{thm_CLT_mf} can be obtained by retracing the computations performed for proving Theorems~\ref{thm_LLN} and~\ref{thm_CLT}.
Indeed, both strategies we used to prove the latter results rely on the knowledge of the limiting free energy (together with its phase diagram) and, thanks to Theorem \ref{fe_meanf_exp}, this information is precisely the same for the mean-field model \eqref{H_mf}
and the edge-triangle model \eqref{Hamilt_ERG}.

Compared with the one of the analogous Theorem~\ref{Thm_convergence_in_distribution}, the proof of Theorem \ref{Thm_convergence_in_distribution_mf} follows by direct computation. In addition, in the present simplified setting we are able to provide the coefficients of the mixture of delta measures describing the infinite volume distribution of the edge density. We give the proof below.

\begin{proof}[Proof of Theorem~\ref{Thm_convergence_in_distribution_mf}.]
We will determine the limit of $\mathbb{\bar{E}}_{n;\alpha,h} \left[ \varphi \left( 2S_n/n^2 \right)\right]$ for any continuous and bounded real function $\varphi$. Recall \eqref{Gibbs_m} and notice that
$$
\bar{\E}_{n;\a,h}\left[\varphi \left(\frac{2S_n}{n^2}\right) \right]
= \sum_{m\in\G_n} \varphi(m) \, \frac{\cN_m e^{n^2 (\frac{\alpha}{6}m^{3} + \frac{h}{2}m)}}{\bar{Z}_{n;\alpha,h}}.
$$
Implementing the same procedure as in the proof of Lemma~\ref{lemma_Z},
and adopting the same notation therein, we split this average over the three sets  $B_{u_1^*}$, $B_{u_2^*}$ and $\cC$, and analyze them separately.
\vspace{0.3cm}

\noindent\emph{Average over the neighborhood of $u_1^*$.}
From the Stirling's approximation \eqref{Stirling} and the identity \eqref{claim},
we get
\be\label{uffa}
\begin{split}
\sum_{m\in B_{u_1^*}} \varphi(m) \, \frac{\cN_m e^{n^2 (\frac{\alpha}{6}m^{3} + \frac{h}{2}m)}}{\bar{Z}_{n;\a,h}}
&=\sum_{m\in B_{u_1^*}} \frac{\varphi(m)}{\sqrt{m(1-m)}}
\frac{e^{-n^2 \left( f_{\a,h} - g_{\a,h}(m)\right)}}{D_1^{(n)}+D_2^{(n)}} (1+o(1))\,.
\end{split}
\ee
Notice that by the definitions of $D_1^{(n)}$ and $D_2^{(n)}$,
and considering the explicit value of $c_1$ and $c_2$, given by
$c_i= -\frac{g_{\a,h}''(u_i^*)}{2} = \frac{1-2\a(u_i^*)^2(1-u_i^*)}{4u_i^*(1-u_i^*)}$
for $i=1,2$, it holds that
\be\label{D12}
\begin{split}
& D_1^{(n)} \, \xrightarrow{\;\; n \to +\infty \;\;} \,
D_1:= \frac{1}{\sqrt{u_1^* \left(1-u_1^* \right)}} \int_{\R} e^{- c_1 x^2} dx
=  2\sqrt{\pi\left({1-2\alpha \left( u_1^*\right)^2(1-u_1^*)}\right)^{-1}}\\
& D_2^{(n)} \, \xrightarrow{\;\; n \to +\infty \;\;} \,
D_2:=\frac{1}{\sqrt{u_2^* \left(1-u_2^* \right)}} \int_{\R} e^{- c_2 x^2} dx
=2 \sqrt{\pi\left({1-2\alpha \left( u_2^*\right)^2(1-u_2^*)}\right)^{-1}} \,.
\end{split}	
\ee
Hence, by the Taylor's expansion  \eqref{taylor_g_M}
and  the change of variable $x=n(m-u_1^*)$,
we can conclude
\be\nonumber
\begin{split}
\sum_{m\in B_{u_1^*}} \varphi(m) \, \frac{\cN_m e^{n^2 (\frac{\alpha}{6}m^{3} + \frac{h}{2}m)}}{\bar{Z}_{n;\a,h}}
&=
\sum_{x\in R^{(n)}} \frac{\varphi(u_1^* +\frac{x}{n})}{\sqrt{(u_1^* +\frac{x}{n})(1-u_1^* -\frac{x}{n})}}
\frac{e^{-c_1 x^2 + \frac{k_1}{n} x^3}}{D_1^{(n)}+D_2^{(n)}} (1+o(1))\\
&\xrightarrow{\;\; n \to +\infty \;\;} \,
\varphi(u_1^*) \frac{D_1}{D_1+ D_2}\,.
\end{split}
\ee
\vspace{0.3cm}

\noindent\emph{Average over the neighborhood of $u_2^*$.}
By the same arguments as for the average over the neighborhood of $u_1^*$, we get
\be\nonumber
\begin{split}
\sum_{m\in B_{u_2^*}} \varphi(m) \, \frac{\cN_m e^{n^2 (\frac{\alpha}{6}m^{3} + \frac{h}{2}m)}}{\bar{Z}_{n;\a,h}}
&=
\sum_{x\in R^{(n)}} \frac{\varphi(u_2^* +\frac{x}{n})}{\sqrt{(u_2^* +\frac{x}{n})(1-u_2^* -\frac{x}{n})}}
\frac{e^{-c_2 x^2 + \frac{k_2}{n} x^3}}{D_1^{(n)}+D_2^{(n)}} (1+o(1))\\
&\xrightarrow{\;\; n \to +\infty \;\;} \,
\varphi(u_2^*) \frac{D_2}{D_1+ D_2}\,.
\end{split}
\ee
\vspace{0.3cm}

\noindent \emph{Average over the complement set $\cC$.}
Since $\varphi$  is a bounded function, from \eqref{somma3}
it follows that the sum over $\cC$ is exponentially small in $n$,
and hence negligible in the large $n$ limit.
\vspace{0.3cm}

Summing up the three contributions, we obtain
$$\lim_{n\to+\infty} \bar{\E}_{n;\a,h}\left[\varphi \left(\frac{2S_n}{n^2}\right) \right]
= \varphi(u_1^*) \frac{D_1}{D_1+ D_2}+\varphi(u_2^*) \frac{D_2}{D_1+ D_2}\,.
$$
\end{proof}

For the mean-field model, we can strengthen the result of Proposition \ref{prop_speed}
as follows.
\begin{proposition}\label{prop_speed_mmf}
For all $(\alpha,h) \in \mathcal{U}^{rs}\setminus \{(\a_c,h_c)\}$,
\begin{equation}\label{speed_mmf1}
\lim_{n\to+\infty}
n \cdot\bar{\mathbb{E}}_{n;\alpha,h}\left( \left\vert \frac{2S_n}{n^2} - u^*(\a,h) \right\vert\right)
=\E(|X|)\,,
\end{equation}
where $X$ is a real Gaussian random variable with Lebesgue density
$\ell(x)\propto e^{-c_0 x^2}$  and
$c_0 
= \frac{1-2\alpha(u^*)^2(1-u^*)}{4u^*(1-u^*)} > 0$.
Moreover, at the critical point, it holds
\begin{equation}\label{speed_mmf2}
\lim_{n\to+\infty}
\sqrt{n}\cdot\bar{\mathbb{E}}_{n;\alpha_c,h_c}
\left( \left\vert \frac{2S_n}{n^2} - u^*(\a_c,h_c) \right\vert\right)
=\E(|Y|)\,,
\end{equation}
where $Y$ is a real random variable with Lebesgue density
$\ell^c(y)\propto e^{-\frac{81}{64}y^4}$.
\end{proposition}

\begin{proof}
Let $(\alpha,h) \in \mathcal{U}^{rs}$ and $u^*=u^*(\a,h)$ be the unique solution
of the scalar problem \eqref{free_energy}.
We keep using the notation introduced in Lemma \ref{lemma_Z}.
Following the same line of argument as in the proof of Theorem
\ref{Thm_convergence_in_distribution_mf}, we obtain
\be\label{prima}
\begin{split}
\bar{\E}_{n;\a,h}\left( \left\vert\frac{2S_n}{n^2}-u^*\right\vert\right)
&= \sum_{m\in B_{u^*}}|m-u^*|
\frac{\cN_m e^{n^{2}(\frac{\alpha}{6}m^{3} + \frac{h}{2}m)}}
{\bar{Z}_{n;\a,h}}(1+o(1))\\
&= \sum_{m\in B_{u^*}}\frac{|m-u^*|}{\sqrt{m(1-m)}}\frac{e^{-n^2 (f_{\a,h}-g_{\a,h}(m))}}
{D^{(n)}(\a,h)}(1+o(1))\,,
\end{split}
\ee
where the second identity is due to Lemma \ref{lemma_Z} and the  Stirling's approximation \eqref{Stirling}. Observe that the infinitesimal correction appearing
in \eqref{prima} is exponentially small in $n$, as it comes from the average
over the complement set $\cC=\G_n\setminus B_{u^*}$ (recall \eqref{somma3})
and Stirling's approximation.

To start, focus on the case $(\alpha,h) \in \mathcal{U}^{rs}\setminus \{(\a_c,h_c)\}$.
By means of the Taylor expansion \eqref{taylor_g} and the change of variable $x= n(m-u^*)$,
 we write
\begin{equation}\nonumber
\begin{split}
n \cdot \bar{\E}_{n;\a,h}\left( \left\vert\frac{2S_n}{n^2}-u^*\right\vert\right)&=
n \cdot \sum_{m\in B_{u^*}}\frac{|m-u^*|}{\sqrt{m(1-m)}}
\frac{e^{-n^{2}c_0(m-u^*)^2 + n^2 k_0(m-u^*)^3}}
{D_*^{(n)}}(1+o(1))\\
&=
\sum_{x\in R^{(n)}}
|x|
\cdot \frac{e^{-c_0x^2 + \frac{k_0}{n}x^3 }}
{\sqrt{(u^*+\frac{x}{n})(1-u^*-\frac{x}{n})}\cdot D_*^{(n)}}(1+o(1))\,,
\end{split}
\end{equation}
where $R^{(n)}:=\left\{ -n^{1-\d}, -n^{1-\d}+\frac{2}{n},\ldots, n^{1-\d} \right\}$.

At last, let  $(X_n)_{n \in \mathbb{N}}$ be a sequence of  random variables
with probability density
$$
\ell_n(x):=
\frac{e^{-c_0x^2 + \frac{k_0}{n}x^3 }}
{{\sqrt{(u^*+\frac{x}{n})(1-u^*-\frac{x}{n})}}\cdot D_*^{(n)}}\1_{R^{(n)}} (x)\,,\quad x\in\R\,,
$$
so that, from the last display, we have
$n\cdot\bar{\E}_{n;\a,h}\left( \left\vert\frac{2S_n}{n^2}-u^*\right\vert\right)= \E(|X_n|)(1+o(1))$.
Notice that, as $n\to+\infty$,  $X_n \xrightarrow{\;\; d \;\;} X$, where $X$ 
is a real Gaussian random variable with density
$$
\ell(x)=  \sqrt{\frac{c_0}{\pi}}e^{-c_0x^2}\,,\quad x\in\R\,.
$$
Moreover, the random variables $X_n$'s have finite exponential  moments for any
sufficiently large $n$.
Therefore, by dominated convergence, we obtain
$$n\cdot \bar{\E}_{n;\a,h}\left( \left\vert\frac{2S_n}{n^2}-u^*\right\vert\right)
= \E(|X_n|)(1+o(1)) \, \xrightarrow{\;\; n \to +\infty \;\;} \, \E(|X|)\,.$$

The proof of statement \eqref{speed_mmf2}, relative to the critical point
$(\a_c,h_c)$, runs similarly; the major difference being that
  the Taylor's approximation \eqref{taylor_g} brings the term $(m-u^*)^4$
  in the exponent of \eqref{prima}.
After the change of variable $y=\sqrt{n}(m-u^*)$, one gets

\begin{equation}\nonumber
\begin{split}
\sqrt{n}\cdot\bar{\E}_{n;\alpha_c,h_c}
\left( \left\vert \frac{2S_n}{n^2} - u^* \right\vert\right)
&=
\sum_{y\in R_c^{(n)}}
|y|
\cdot
\frac{e^{-\frac{81}{64}y^4 + \frac{k_c}{\sqrt n} y^5}}
{\sqrt{(u^*+\frac{y}{\sqrt n})(1-u^*-\frac{y}{\sqrt n})}\cdot  D_c^{(n)}} (1+o(1))
\\
&
=\E(|Y_n|)(1+o(1))\,,
\end{split}
\end{equation}
where $R_c^{(n)}:= \left\{-n^{1/2-\delta}, -n^{1/2-\delta}+\frac{2}{n^{3/2}}, \dots, n^{1/2-\delta}\right\}$ and $Y_n$ is a real random variable with Lebesgue
density

\be\label{density_yn}
\ell^c_n(y):=
\frac{e^{-\frac{81}{64}y^4  + \frac{k_c}{\sqrt{n}}y^5 }}
{\sqrt{(u^*+\frac{y}{\sqrt n})(1-u^*-\frac{y}{\sqrt n})} \cdot D_c^{(n)}}\1_{R_c^{(n)}} (y)\,,\quad y\in\R\,.
\ee
The random variable $Y_n$ has finite  exponential moments for any
sufficiently large $n$.
Since, as $n\to+\infty$,   $Y_n \xrightarrow{\;\; d \;\;} Y$,
 where $Y$ is a real random variable with symmetric density
 $\ell^c(y)\propto e^{-\frac{81}{64}y^4}$,
the statement \eqref{speed_mf2}  follows by dominated convergence:
$$\sqrt{n}\cdot \bar{\E}_{n;\alpha_c,h_c}
\left( \left\vert \frac{2S_n}{n^2} - u^* \right\vert\right)
= \E(|Y_n|)(1+o(1)) \, \xrightarrow{\;\; n \to +\infty \;\;} \, \E(|Y|)\,.$$

\end{proof}

The strategy of the proof of Proposition \ref{prop_speed_mmf}
can be implemented to prove the following.
\begin{corollary}\label{cor_speed_mf}
For all $(\alpha,h) \in \mathcal{U}^{rs}\setminus \{(\a_c,h_c)\}$, we have
\begin{equation}\label{speed_mf1}
\lim_{n\to+\infty} n \cdot (\bar m_n(\a,h)-u^*(\a,h)) =0\,,
\end{equation}
while at the critical point it holds
\begin{equation}\label{speed_mf2}
\lim_{n\to+\infty} \sqrt{n}\cdot(\bar m_n(\a_c,h_c)-u^*(\a_c,h_c)) =0\,.
\end{equation}
\end{corollary}

\begin{proof}
Recall that $\bar m_n(\a,h)= \bar{\E}_{n;\a,h}\left( \frac{2S_n}{n^2}\right)$.
Following step by step the proof of Proposition \ref{prop_speed_mmf},
with the same notation introduced there, we  get that:
\begin{itemize}
\item
 for all
$(\alpha,h) \in \mathcal{U}^{rs}\setminus \{(\a_c,h_c)\}$, then
$$n\cdot\left( \bar m_n(\a,h)-u^*(\a,h)\right)
= \E(X_n)(1+o(1)) \, \xrightarrow{\;\; n \to +\infty \;\;} \, \E(X)=0,$$
\item
if $(\alpha,h)=(\a_c,h_c)$, then
$$\sqrt{n}\cdot\left( \bar m_n(\a_c,h_c)-u^*(\a_c,h_c)\right)
= \E(Y_n)(1+o(1)) \, \xrightarrow{\;\; n \to +\infty \;\;} \, \E(Y)=0,$$
\end{itemize}
that ends the proof.
\end{proof}

Notice that the central limit theorem stated in Theorem~\ref{thm_CLT_mf} is true even if we replace $\bar{m}_n(\alpha,h)$ by $u^*(\alpha,h)$. Indeed, in view of \eqref{CLT_mf} and \eqref{speed_mf1}, the decomposition
\[
\sqrt{2} \, \frac{S_n-\frac{n^2}{2}u^*(\alpha,h)}{n} = \sqrt{2} \, \frac{S_n-\frac{n^2}{2}\bar{m}_n(\alpha,h)}{n} + \frac{n}{\sqrt{2}} (\bar{m}_n(\alpha,h)-u^*(\alpha,h))
\]
and Slutsky's theorem prove the assertion.

Now we provide a non-standard central limit theorem at the critical point.

\begin{theorem}
[Non-standard CLT for $S_n$ w.r.t. $\bar{\mathbb{P}}_{n;\alpha_c,h_c}$]\label{Thm_non-standard_CLT}
If $(\a,h)=(\alpha_c,h_c)$, it holds
\[	
2 \, \frac{S_n - \frac{n^{2}}{2}\bar m_n(\a_c,h_c)}{n^{3/2}} \xrightarrow{\;\;\mathrm{d}\;\;}{} Y \quad \text{ w.r.t. } \bar{\mathbb{P}}_{n;\alpha_c,h_c},\text{ as } n \to +\infty,
\]
where $Y$ is a real random variable with Lebesgue density $\ell^c(y)\propto e^{-\frac{81}{64}y^{4}}$.	
\end{theorem}

\proof
Recall that $\alpha_{c}=\frac{27}{8}$, $h_{c}= \ln 2 -\frac{3}{2}$
and $u^{*}(\alpha_c,h_{c})=\frac{2}{3}$.
We keep using the notation introduced in Lemma \ref{lemma_Z} and, not to
clutter too much our formulas,  we  write $u^{*}$ in place of
$u^{*}(\alpha_c,h_{c})$ along the computations.

Consider the trivial decomposition
\be\label{equivalenza}
2 \, \frac{S_n - \frac{n^{2}}{2}\bar m_n(\a_c,h_c)}{n^{3/2}}
= 2 \, \frac{S_n - \frac{n^{2}}{2}u^{*}}{n^{3/2}}
+ \sqrt{n}\left(u^{*} - \bar{m}_n(\a_c,h_c)\right).
\ee

Due to Corollary \ref{cor_speed_mf},
it is enough to study the convergence in distribution of
$2(S_n - \frac{n^{2}}{2}u^*)/n^{3/2}$.
To this purpose, we analyze the moment generating function
of this random variable with the final goal of showing that, for any $t\in\R$,
\be\label{goal}
\bar{M}_n (t):=\bar{\mathbb{E}}_{n;\alpha_c,h_c}\left(e^{t\left[2\frac{S_{n}-\frac{n^{2}}{2}u^{*}}{n^{3/2}}\right]}\right) \, \xrightarrow{\;\; n \to +\infty \;\;} \,
\int_{\mathbb{R}}e^{ty}\ell^c(y)dy,
\ee
where $\ell^c$ is the density given in the statement. The proof follows the ideas developed in the proof of Proposition \ref{prop_speed_mmf}.
Observe that
\be\label{base}
\begin{split}
\bar{M}_n (t) &= \sum_{m\in B_{u^*}}
\frac{\cN_m e^{ t\sqrt{n}(m-u^{*})+ n^{2}(\frac{\alpha_c}{6}m^{3} + \frac{h_c}{2}m)}}
{\bar{Z}_{n;\a_c,h_c}}(1+o(1))\\
&= \sum_{m\in B_{u^*}}\frac{1}{\sqrt{m(1-m)}}\frac{e^{t\sqrt{n}(m-u^{*})-n^2 (f_{\a_c,h_c}-g_{\a_c,h_c}(m))}}
{D_c^{(n)}}(1+o(1))\,,
\end{split}
\ee
where the second identity follows from Lemma \ref{lemma_Z} and the  Stirling's approximation \eqref{Stirling}.
We mention that the definition of  $B_{u^*}$ as
the $n^{-\d}$-neighborhood of $u^*$, together with the hypothesis $\d< \frac{3}{8}$,
are crucial to ensure that the contribution of the average over
the complement set $\cC= \G_n\setminus B_{u^*}$  is negligible
for large $n$.
Indeed, from  the Stirling's approximation \eqref{Stirling-rough}
and the consequent  rough bound $\bar{Z}_{n;\a,h} \geq cn^{-1} e^{n^2 f_{\a,h}}$,
we get
\begin{align*}
\sum_{m\in \cC}
\frac{\cN_m e^{t\sqrt{n}(m-u^{*}) +n^{2}(\frac{\alpha_c}{6}m^{3} + \frac{h_c}{2}m)}}
{\bar{Z}_{n;\a_c,h_c}}
&\leq C c^{-1} n^4 e^{t\sqrt{n} -n^2 (f_{\a_c,h_c}- \max_{m\in\cC} g_{\a_c,h_c}(m))}\\
&\leq e^{- k n^{-4\d} },
\end{align*}
where the last inequality holds for some $k>0$ and any  large enough $n$.
We used that $|m-u^*|^4>n^{-4\d}$ for $m\in\cC$, and that
$2-4\d>1/2$, whenever $\d<3/8$.

Inserting the Taylor's expansion \eqref{taylor_g} in \eqref{prima},
and making the change of variable $y=(m-u^{*})\sqrt{n}$,
we  get
\begin{align}\nonumber
\bar{M}_n (t)&=
\sum_{m\in B_{u^*}}\frac{1}{\sqrt{m(1-m)}}
\frac{e^{t\sqrt{n}(m-u^{*})-\frac{81}{64}n^{2}(m-u^{*})^{4}+ k_c n^{2}(m-u^{*})^{5}}}
{D_c^{(n)}}(1+o(1))\\
&\nonumber =
\sum_{y\in R_c^{(n)}}
e^{ty}\frac{e^{ -\frac{81}{64}y^{4} +  k_c\frac{y^5}{\sqrt{n}}}}
{ \sqrt{(u^*+\frac{y}{\sqrt{n}})(1-u^*-\frac{y}{\sqrt{n}})} D_c^{(n)}}
(1+o(1))\,.
\end{align}
In the last line we can recognize the probability density
$\ell_n^c$ of the random variable $Y_n$ we introduced in \eqref{density_yn}.
Thus we can conclude that
$\bar{M}_n (t)= \E\left(e^{tY_n}\right)(1+o(1))\,.$
Since $Y_n \xrightarrow{\;\; d \;\;} Y$, and $Y_n$ and $Y$
have finite exponential moments, convergence \eqref{goal} follows for all $t\in\R$.
\endproof

\begin{remark}\label{rmk:variance_mfm}
The very same strategy adopted to prove Theorem~\ref{Thm_non-standard_CLT} can be used to give an alternative proof of the central limit theorem stated in Theorem~\ref{thm_CLT_mf}. 
The approach of characterizing the moment generating function has the advantage of allowing to obtain the \emph{explicit} value of the variance \eqref{variance_CLT_mfm}. In particular, it yields
\begin{equation}\label{explicit_variance_CLT_mfm}
v(\alpha,h) = \left(4c_0\right)^{-1} = \frac{u^*\left(1-u^*\right)}{1-2\alpha \left(u^*\right)^2\left(1-u^*\right)}.
\end{equation}
In view of Theorem~\ref{fe_meanf_exp}, the edge-triangle model and its mean-field approximation have the same infinite volume free energy, and therefore the same maximizer $u^*(\alpha,h)$, for $(\alpha,h) \in \mathcal{U}^{rs}$. 
As a consequence, thanks to \eqref{explicit_variance_CLT_mfm}, we obtain the value of the variance \eqref{variance_CLT_mfm} in the central limit theorem for the edge-triangle model. 
\end{remark}

We conclude the study of the mean-field model by analyzing further the behavior of the system on the critical curve.
In particular, we want to analyze the limiting behavior of the edge density when confined in a suitable neighborhood of a solution of the scalar problem \eqref{free_energy}.
To this purpose, we introduce two conditional measures.\\

Consider $(\alpha,h) \in \mathcal{M}^{rs}$ and let $u_i^*=u_i^*(\alpha,h)$ ($i=1,2$) be the solutions of the scalar problem \eqref{free_energy}.
Moreover, for any $\delta \in (0,1)$, let $B_{u_i^*}$ ($i=1,2$) be the neighborhoods introduced in \eqref{neighborhoods}. 
For $i=1,2$ we define the conditional probability measures
\begin{equation}\label{conditional_measure_mfm}
\hat{\mathbb{P}}_{n;\alpha,h}^{(i)}\left(\frac{2S_n}{n^2} \in A \right) := \bar{\mathbb{P}}_{n;\alpha,h} \left(\frac{2S_n}{n^2} \in A \left\vert  \frac{2S_n}{n^2} \in B_{u_i^*} \right.\right), \qquad \text{ for } A \subseteq \Gamma_n,
\end{equation}
and we denote the corresponding averages by $\hat{\mathbb{E}}_{n;\alpha,h}^{(i)}$. 

Recall that, by definition, the neighborhoods $B_{u_i^*}$ depend on $\d$, and so do the conditional measures $\hat{\mathbb{P}}_{n;\alpha,h}^{(i)}$.
Not to clutter our formulas, we prefer to leave this dependence implicit. We remark that all the following results are valid for any choice of $\delta \in (0,1)$.

We want to prove standard limit theorems for the edge density under the probability measures $\hat{\mathbb{P}}_{n;\alpha,h}^{(i)}$ ($i=1,2$).
We start by proving two auxiliary results concerning the speed of convergence of the edge density to its average value.

\begin{proposition}\label{prop_speed_conditional_mfm}
For $i=1,2$ and for all $(\alpha,h) \in \mathcal{M}^{rs}$,
\begin{equation}\label{speed_1_conditional_mfm}
\lim_{n \to +\infty} n \cdot \hat{\mathbb{E}}_{n;\alpha,h}^{(i)} \left( \left\vert \frac{2S_n}{n^2} - u_i^*(\alpha,h) \right\vert \right) = \mathbb{E}\left(\left\vert X^{(i)} \right\vert\right),
\end{equation}
where $X^{(i)}$ is a real Gaussian random variable with Lebesgue density $\ell_i(x) \propto e^{-c_i x^2}$ and 
$c_i
= \frac{1-2\alpha(u_i^*)^2(1-u_i^*)}{4u_i^*(1-u_i^*)}>0$.
\end{proposition}

\proof
For the sake of concreteness, we stick with the case $i=1$, the other being similar. 
Let $(\alpha,h) \in \mathcal{M}^{rs}$ and $u_1^*=u_1^*(\a,h)$ be a solution of the scalar problem \eqref{free_energy}.
We keep using the notation introduced in Lemma \ref{lemma_Z} and we mimic the proof of Proposition~\ref{prop_speed_mmf}.
We have
\begin{align*}
\hat{\mathbb{E}}_{n;\a,h}^{(1)}\left( \left\vert \frac{2S_n}{n^2} - u_1^* \right\vert \right) 
&= \sum_{m\in B_{u^*}} \left\vert m-u_1^* \right\vert \frac{\cN_m \, e^{n^{2}(\frac{\alpha}{6}m^{3} + \frac{h}{2}m)}} {\bar{Z}_{n;\a,h}(B_{u_1^*})}\\
&= \sum_{m\in B_{u^*}}\frac{|m-u^*|}{\sqrt{m(1-m)}}\frac{e^{-n^2 (f_{\a,h}-g_{\a,h}(m))}}
{D_1^{(n)}}(1+o(1)),
\end{align*}
where the identities are due to \eqref{somma1} and the  Stirling's approximation \eqref{Stirling}. 
Observe that the infinitesimal correction in the last line of the previous display is exponentially small in $n$, as it comes from the Stirling's approximation.

By means of the Taylor expansion \eqref{taylor_g_M} and the change of variable $x= n(m-u^*)$, we obtain
\[
n \cdot \hat{\mathbb{E}}^{(1)}_{n;\a,h}\left( \left\vert\frac{2S_n}{n^2}-u_1^*\right\vert\right) =
\sum_{x\in R^{(n)}}
|x|
\cdot \frac{e^{-c_1 x^2 + \frac{k_1}{n}x^3 }}
{\sqrt{(u^*+\frac{x}{n})(1-u^*-\frac{x}{n})}\cdot D_1^{(n)}}(1+o(1)),
\]
where $R^{(n)}:=\left\{ -n^{1-\d}, -n^{1-\d}+\frac{2}{n},\ldots, n^{1-\d} \right\}$.
In other words, if $(X_n)_{n \in \mathbb{N}}$ is a sequence of random variables
with probability density
$$
\ell_n^{(1)}(x):=
\frac{e^{-c_1 x^2 + \frac{k_1}{n} x^3 }}
{{\sqrt{(u^*+\frac{x}{n})(1-u^*-\frac{x}{n})}}\cdot D_1^{(n)}} \, \1_{R^{(n)}} (x),\quad x\in\R,
$$
from the last display, we can write
$n \cdot \hat{\mathbb{E}}^{(1)}_{n;\a,h}\left( \left\vert\frac{2S_n}{n^2}-u_1^*\right\vert\right) = \E(|X_n|)(1+o(1))$.
Notice that, as $n\to+\infty$,  $X_n \xrightarrow{\;\; d \;\;} X^{(1)}$, where $X^{(1)}$ is a real Gaussian random variable with density
$$
\ell_1(x)=  \sqrt{\frac{c_1}{\pi}}e^{-c_1 x^2},\quad x\in\R.
$$
Moreover, as the random variables $X_n$'s have finite exponential moments for any
sufficiently large $n$, by dominated convergence, we get
$$n \cdot \hat{\mathbb{E}}^{(1)}_{n;\a,h}\left( \left\vert\frac{2S_n}{n^2}-u_1^*\right\vert\right)
= \E(|X_n|)(1+o(1)) \, \xrightarrow{\;\; n \to +\infty \;\;} \, \E\left(\left\vert X^{(1)}\right\vert\right).$$
\endproof

\begin{corollary}\label{cor_speed_conditional_mfm}
Set $\hat{m}_n^{(i)}(\alpha,h):=\hat{\mathbb{E}}^{(i)}_{n;\alpha,h}\left(\frac{2S_n}{n^2}\right)$ ($i=1,2$). For $i=1,2$ and for all $(\alpha,h) \in \mathcal{M}^{rs}$, we have
\[
\lim_{n \to +\infty} n \cdot \left(\hat{m}_n^{(i)}(\alpha,h)-u_i^*(\alpha,h)\right) = 0.
\]
\end{corollary}

\proof
Following step by step the proof of Proposition \ref{prop_speed_conditional_mfm},
with the same notation therein, we  get that, for all $(\alpha,h) \in \mathcal{M}^{rs}$, 
\[
n \cdot \left( \hat{m}_n^{(1)}(\a,h)-u_1^*(\a,h)\right)
= \E(X_n)(1+o(1)) \, \xrightarrow{\;\; n \to +\infty \;\;} \, \E \left( X^{(1)} \right)=0.
\]
\endproof

The next theorem is the analog of Theorems~\ref{thm_LLN_mf} and~\ref{thm_CLT_mf}, but it is obtained under the constraint (\emph{conditioning}) that the edge density is close to one of the maximizers of the scalar problem \eqref{free_energy}.

\begin{theorem}[Conditional LLN and CLT]
\label{thm_conditional_limit_theorems_mfm}
For $i=1,2$ and for all $(\alpha,h) \in \mathcal{M}^{rs}$, it holds
\begin{equation}\label{conditional_lln_mfm}
\frac{2S_n}{n^2} \, \xrightarrow{\;\;\mathrm{a.s.}\;\;}{} \, u_i^{*}(\alpha,h)   \quad \text{ w.r.t. } \hat{\mathbb{P}}^{(i)}_{\alpha,h}, \text{ as } n \to +\infty,
\end{equation}
and
\begin{equation}\label{conditional_clt_mfm}
\sqrt{2} \, \frac{S_n - \frac{n^2}{2} \hat{m}_n^{(i)}(\alpha,h)}{n} \, \xrightarrow{\;\;\mathrm{d}\;\;}{} \, \mathcal{N}(0,v_i(\alpha,h))   \quad \text{ w.r.t. } \hat{\mathbb{P}}^{(i)}_{n;\alpha,h}, \text{ as } n \to +\infty,
\end{equation}
where $\mathcal{N}(0,v_i(\alpha,h))$ is a centered Gaussian distribution with variance
\[
v_i(\alpha,h) = \frac{u_i^*(\alpha,h)[1-u_i^*(\alpha,h)]}{1-2\alpha[u_i^*(\alpha,h)]^2[1-u_i^*(\alpha,h)]}.
\]
\end{theorem}

\proof
We prove only the $i=1$ case. 
In the sequel, not to clutter too much our formulas, we will drop the dependence on $\alpha$ and $h$ in the function $u_1^*=u_1^*(\alpha,h)$.

Firstly, we show the statement \eqref{conditional_lln_mfm}.
%
%
We prove that, under $\hat{\mathbb{P}}^{(1)}_{n;\alpha,h}$, the sequence $(2S_n/n^2)_{n \geq 1}$ converges exponentially to $u_1^*$, according to Definition~\ref{def:exp_conv}.
The almost sure convergence then follows by a standard Borel-Cantelli argument (see~\cite{E}, Thm.~II.6.4).

Let $\mathcal{R} \equiv \mathcal{R}(\eta;n) := \left\{ m \in \Gamma_n: \eta \leq \left\vert m - u_1^* \right\vert < n^{-\delta} \right\}$ be the circular crown with minor radius $\eta$ and major radius $n^{-\delta}$. 
For any $\eta > 0$, we have
\begin{align*}
\hat{\mathbb{P}}^{(1)}_{n;\alpha,h} \left( \left\vert \frac{2S_n}{n^2} - u_1^*\right\vert \geq \eta \right) &\leq \sum_{m \in \mathcal{R}}  \frac{\mathcal{N}_m e^{n^2 \left(\frac{\alpha}{6}m^3+\frac{h}{2}m\right)}}{\bar{Z}_{n;\alpha,h}(B_{u_1^*})} \\
&\leq  C c^{-1} n^4 e^{-n^2 (f_{\alpha,h}-\max_{m \in \mathcal{R}} g_{\alpha,h}(m))}\\
&\leq C c^{-1} n^4 e^{-n^2 \min_{m \in \mathcal{R}} (f_{\alpha,h} - g_{\alpha,h}(m))},
\end{align*}
where the second-last inequality is due to the Stirling's approximation \eqref{Stirling-rough} and the resulting estimate $\bar{Z}_{n;\alpha,h}(B_{u_1^*}) \geq cn^{-1} e^{n^2 f_{\alpha,h}}$.
For sufficiently large $n$, the function $f_{\alpha,h} - g_{\alpha,h}(m)$ restricted to the neighborhood $B_{u_1^*}$ is positive, convex and admits a unique zero at $m=u_1^*$ (see~\cite{RY}, Prop.~3.2). 
As a consequence, we obtain $\min_{m \in \mathcal{R}} (f_{\alpha,h} - g_{\alpha,h}(m)) = \min\{f_{\alpha,h} - g_{\alpha,h}(u_1^*-\eta), f_{\alpha,h} - g_{\alpha,h}(u_1^*+\eta)\}>0$, giving the conclusion.

We turn now to the convergence \eqref{conditional_clt_mfm}. 
The proof of this statement runs similarly to the proof of Theorem~\ref{Thm_non-standard_CLT}.
Consider the trivial decomposition
$$\sqrt{2} \, \frac{S_n - \frac{n^{2}}{2}\hat{m}_n^{(1)}(\a,h)}{n}
= \frac{n}{\sqrt{2}} \left( \frac{2S_n}{n^2} - u_1^{*}\right) 
+ \frac{n}{\sqrt{2}} \left(u_1^{*} - \hat{m}^{(1)}_n(\a,h)\right).$$
By Corollary \ref{cor_speed_conditional_mfm}, it suffices to study the convergence in distribution of the random variable $\frac{n}{\sqrt{2}} \left( \frac{2S_n}{n^2} - u_1^{*}\right)$.
To this purpose, we analyze its moment generating function 
\[
\hat{M}_n(t) := \hat{\mathbb{E}}^{(1)}_{n;\alpha,h}\left(e^{\frac{tn}{\sqrt{2}} \left( \frac{2S_n}{n^2} - u_1^{*}\right)}\right),
\]
and we show that, for any $t\in\R$, it converges to the moment generating function of a centered Gaussian having variance $v_1(\alpha,h)$. 

Observe that
\begin{align*}
\hat{M}_n (t) &= \sum_{m\in B_{u_1^*}}
\frac{\cN_m e^{ \frac{tn}{\sqrt{2}}(m-u_1^{*})+ n^{2}(\frac{\alpha}{6}m^{3} + \frac{h}{2}m)}}
{\bar{Z}_{n;\a,h}(B_{u_1^*})}\\
&= \sum_{m\in B_{u_1^*}}\frac{1}{\sqrt{m(1-m)}}\frac{e^{\frac{tn}{\sqrt{2}}(m-u_1^{*})-n^2 (f_{\a,h}-g_{\a,h}(m))}}
{D_1^{(n)}}(1+o(1)),
\end{align*}
where the second identity follows from the identity \eqref{somma1} and the  Stirling's approximation \eqref{Stirling}.
By using the Taylor's expansion \eqref{taylor_g_M} and making the change of variable $z=\frac{n}{\sqrt{2}}(m-u_1^{*})$, we  get
\begin{align*}
\hat{M}_n (t) &=
\sum_{z \in \tilde{R}^{(n)}}
e^{tz} \, \frac{e^{ -2 c_1 z^{2} +  2 \sqrt{2} k_1 \frac{z^3}{n}}}
{\sqrt{\left(u^*+\frac{\sqrt{2}z}{n}\right)\left(1-u^*-\frac{\sqrt{2}z}{n}\right)} \cdot D_1^{(n)}}
(1+o(1))\\[.2cm]
&= \mathbb{E}\left( e^{tZ_n} \right)(1+o(1)),
\end{align*}
where $\tilde{R}^{(n)} := \left\{ -\frac{n^{1-\delta}}{\sqrt{2}}, -\frac{n^{1-\delta}}{\sqrt{2}} + \frac{\sqrt{2}}{n}, \dots, \frac{n^{1-\delta}}{\sqrt{2}}, -\frac{n^{1-\delta}}{\sqrt{2}} \right\}$  and $Z_n$ is a real random variable with Lebesgue density
\[
\tilde{\ell}_n(z) :=
\frac{e^{ -2 c_1 z^{2} +  2 \sqrt{2} k_1 \frac{z^3}{n}}}
{\sqrt{\left(u^*+\frac{\sqrt{2}z}{n}\right)\left(1-u^*-\frac{\sqrt{2}z}{n}\right)} \cdot  D_1^{(n)}} \, \1_{\tilde{R}^{(n)}} (z),\quad z\in\R.
\]
Notice that, as $n \to +\infty$, $Z_n \xrightarrow{\;\; d \;\;} Z$, where $Z$ is a real Gaussian random variable with density
\[
\tilde{\ell}(z) = \sqrt{\frac{2c_1}{\pi}}e^{-2c_1 x^2},\quad x\in\R.
\]
Since $Z$ has finite exponential moments and $Z_n$ as well, at least for sufficiently large $n$, the convergence
\[
\hat{M}_n(t) \xrightarrow{\;\; n \to +\infty \;\;} \int_{\mathbb{R}} e^{tz} \tilde{\ell}(z) \, dz
\]
follows for all $t\in\R$.

\endproof

\subsection{Comparison among the edge-triangle model and its mean-field approximation}
\label{sect:comparison_ET&mfm}

The analysis of the mean-field model carried out in this section has highlighted various points of similarity  among the asymptotic behavior of the edge-triangle model and its mean-field approximation.
These are essentially derived as consequences of the equality between the infinite volume free energies of the two models, stated in the Theorem~\ref{fe_meanf_exp}.  
Thanks to the regularity of $f_{\a,h}$ for parameters in the region of $\mathcal{U}^{rs}\setminus\{(\a_c,h_c)\}$ and some classical statistical mechanics tools, it was possible to derive analogous SLLN and CLT for the edge density in the two models. In turn, this common behavior has allowed us to infer some hidden properties of the edge-triangle
model. For example, by comparison, we have been able to provide an explicit formula for the variance $v(\a,h)$ appearing in Theorem~\ref{thm_CLT} (see Remark~\ref{rmk:variance_mfm}).

The main differences between these models emerge while analyzing the behavior
at criticality: at the critical point $(\a_c,h_c)$ and in the multiplicity region $\mathcal{M}^{rs}$. For parameters in $\mathcal{M}^{rs} \cup \{(\a_c,h_c)\}$, the loss of regularity of $f_{\a,h}$ must be compensated with a deeper probabilistic analysis of the model,  and large deviations estimates, mainly connected to the analysis given in~\cite{CD}, come into play.

However, the mean-field approximation is a simpler model to deal with and we can obtain much more detailed information on its asymptotic behavior.
In particular, the very precise characterization of the partition function given in Lemma~\ref{lemma_Z} goes much beyond a 
large deviation principle for the sequence of measures $(\bar{\mathbb{P}}_{n;\a,h})_{n \geq 1}$,
and it is indeed a key tool for the analysis of the mean-field model at criticality.

Beyond technical difficulties, it is natural to ask whether, at criticality, analogous results hold 
for the edge density of the edge-triangle model and its mean-field approximation.
As a first attempt of comparison, let us focus on the non-standard CLT stated in Theorem~\ref{Thm_non-standard_CLT}, valid
at the critical point $(\a_c,h_c)$, and on its proof. 
Since the result is proved by showing the convergence of the moment generating function given in \eqref{goal}, 
one can try to recover a similar convergence for the edge-triangle model by analyzing the corresponding moment generating function. 
If we set
$$M_n(t):=\mathbb{E}_{n;\alpha_{c},h_{c}}\left(e^{t\left[2\frac{S_{n}
-\frac{n^{2}}{2}u^*}{n^{3/2}}\right]}\right)\,,\quad \forall t\in\R\,,
$$
we would like to prove that $\lim_{n\to\infty}\left(M_n(t)-\bar{M}_n(t)\right)=0$
or, equivalently, that
\begin{equation}\label{fin_goal}
\lim_{n\to\infty} \left(\log M_n(t)-\log \bar{M}_n(t)\right)=0\,.
\end{equation}
Elaborating on the characterization of the cumulant generating function of $S_n$ w.r.t. $\P_{n;\a_c,h_c}$ given in \eqref{clt_decomp} and setting $\bar{c}_{n}(t):= 2\left(\bar{f}_{n;\alpha_{c},h_{c}+ t}- \bar{f}_{n;\alpha_{c},h_{c}}\right)$, the convergence \eqref{fin_goal} 
can be proved by studying the limit, as $n \to\infty$,  of
\be\label{goal_critClt}
\frac{n^2}{2}c_n(t_n) - \frac{n^2}{2}\bar{c}_n(t_n), 
\ee
where $t_n:=2t/n^{3/2}$. 
By expanding the cumulant generating functions around $0$, we obtain that, for some $t^{*}_{n}\in(0,t_n)$,
\begin{equation}\nonumber 
c_{n}(t_n) =2\left(f_{n;\alpha_c;h_c+t_n}- f_{n;\alpha_c;h_c}\right)= 2t_n \partial_h f_{n;\alpha_c;h_c+t^{*}_n}
= t_n m_n(\a_c,h_c+t^{*}_n)
\end{equation}
 and, similarly, $\bar{c}_{n}(t_n)=\frac{t_n}{2} \bar{m}_n(\a_c,h_c+\bar{t}^{*}_n)$. 
 Substituting these expansions into \eqref{goal_critClt}, and inserting the specific value of $t_n$, we then get
\be\label{stima}
\frac{t\sqrt{n}}{2}\left[m_n(\alpha_c,h_c +t^*_n) - \bar{m}_n(\alpha_c,h_c +\bar{t}^*_n)\right] \,.
\ee
Although we do not have direct tools to compare the above difference,
since the average edge densities in the square brackets are averaged 
over distinct measures,
notice that $m_n(\alpha_c,h_c +t^*_n)$ and $\bar{m}_n(\alpha_c,h_c +\bar{t}^*_n)$ share the limit.
Indeed, since for all $n \in \N$ the functions $c_{n}(t)$ and $\bar{c}_{n}(t)$ are convex with finite limit for all $t\in\mathbb{R}$ (see equation \eqref{cgf}), 
and such that $c'_{n}(t^{*}_n)=m_{n}(\alpha_c,h_c+t^{*}_n)$, $\bar{c}'_{n}(\bar{t}^{*}_n)=\bar{m}_{n}(\alpha_c,h_c+\bar{t}^{*}_n)$ and $c'(0)=\bar{c}'(0)=u^{*}$ exist,  
from \cite[Lemma V.7.5]{E} it follows that
$$
\lim_{n\to\infty}m_{n}(\alpha_c,h_c+t^{*}_{n})=
\lim_{n\to\infty}\bar{m}_{n}(\alpha_c,h_c+\bar{t}^{*}_{n})
= \lim_{n\to\infty}m_{n}(\alpha_c,h_c)=u^{*}.
$$
These identities can then be used to rewrite \eqref{stima} as the following sum of four differences
\be\nonumber
\begin{split}
&\frac{t\sqrt{n}}{2}\left[
\left(m_n(\alpha_c,h_c +t^{*}_n) - m_n(\alpha_c,h_c)\right)
+ \left(m_n(\alpha_c,h_c) - u^*\right)
\right]\\
&\quad
- \frac{t\sqrt{n}}{2}\left[
\left(\bar{m}_n(\alpha_c,h_c +\bar{t}^{*}_n) -\bar{m}_n(\alpha_c,h_c)\right)
+\left(\bar{m}_n(\alpha_c,h_c) - u^*\right)
\right]\,.
\end{split}
\ee
But with this rewriting, two different technical problems emerge.
On the one hand, on each of the two lines above, the first difference can be thought of as an approximation of the second derivative in $h$ of the finite volume free energy,
whose limit, as $n\to\infty$, explodes as an unknown function of $n$.
On the other hand, while the second term in the second line 
converges to $0$ due to Corollary \ref{cor_speed_mf}, the convergence of the analogous term in the first line is not guaranteed. 
The technical element missing at this point is a sufficient control of the speed of convergence of $m_{n}(\alpha_c,h_c)$ toward $u^*$,
as apparent by comparing Proposition \ref{prop_speed} with Proposition \ref{prop_speed_mmf}. 
For the same reason, while in the mean-field setting 
$2(S_n - \frac{n^{2}}{2}u^*)/n^{3/2}$ and $2(S_n - \frac{n^{2}}{2}\bar{m}_n(\a_c,h_c))/n^{3/2}$ share the same
limiting law (recall identity \eqref{equivalenza}), 
in the edge-triangle model this might no longer be the case.


\end{document}